\documentclass[11pt, twoside,fleqn]{article}
\usepackage{amsmath,amsthm,amssymb,amsfonts}
\usepackage[dvipsnames]{xcolor}
\usepackage{mathptmx}
\numberwithin{equation}{section}
\newcommand{\rr}{\mathbb{R}}
 \newcommand{\rn}{\mathbb{R}^n}
  \newcommand{\Om}{\Omega}
   \newcommand{\test}{C_c^\infty}
    \newcommand{\prt}{\partial}
     \newcommand{\N}{\mathbb{N}}
      \newcommand{\dx}{\mathrm{d}x}
       \newcommand{\dS}{\mathrm{d}\mathcal{H}^{n-1}}
        \newcommand{\us}{\mathbb{S}^{n-1}}
\newcommand{\al}{\alpha}
 \newcommand{\de}{\delta}
  \newcommand{\e}{\varepsilon}
   \newcommand{\omn}{\omega_{n}}
    \newcommand{\dd}{\textrm{d}}
     \newcommand{\sgn}{\mathrm{sgn}}
\newcommand{\be}{\begin{equation}}
 \newcommand{\ee}{\end{equation}}
  \newcommand{\bea}{\begin{eqnarray}}
   \newcommand{\eea}{\end{eqnarray}}
  \DeclareMathOperator{\Div}{div}
   \DeclareMathOperator{\diam}{diam}

\newtheorem{theorem}{Theorem}[section]
 \newtheorem{proposition}[theorem]{Proposition}
  \newtheorem{lemma}[theorem]{Lemma}
   \newtheorem{corollary}[theorem]{Corollary}
    \theoremstyle{definition}
     \newtheorem{definition}[theorem]{Definition}
      
       \newtheorem{remark}[theorem]{Remark}

\begin{document}
\title{Optimal non-homogeneous improvements\\ for the series expansion of Hardy's inequality}
\author{K. T. Gkikas\footnote{Konstantinos T. Gkikas: Centro de Modelamiento Matem\'{a}tico (UMI 2807 CNRS), Universidad de Chile, Casilla 170 Correo 3, Santiago, Chile; \texttt{kugkikas@gmail.com}}
 \and
  G. Psaradakis\footnote{Georgios Psaradakis: Institut f\"{u}r Mathematik (Lehrst\"{u}hl f\"{u}r Mathematik IV), Universit\"{a}t Mannheim, A5, Mannheim 68131, Deutschland; \texttt{psaradakis@uni-mannheim.de}}}
\date{}
\maketitle

\begin{abstract}

We consider the series expansion of the $L^p$-Hardy inequality of \cite{BFT2}, in the particular case where the distance is taken from an interior point of a bounded domain in $\rn$ and $1<p\neq n$. For $p<n$ we improve it by adding as a remainder term an optimally weighted critical Sobolev norm, generalizing the $p=2$ result of \cite{FT} and settling the open question raised in \cite{BFT1}. For $p>n$ we improve it by adding as a remainder term the optimally weighted H\"{o}lder seminorm, extending the Hardy-Morrey inequality of \cite{Ps} to the series case.

\medskip


\medskip

\end{abstract}
\section{Introduction}
Let $\Om$ be any domain in $\rn$, $n\geq3$, containing the origin. Hardy's inequality asserts that for all $u\in H_0^1(\Om)$
\be\label{Hardy:L^2}
 \int_{\Om}|\nabla u|^2~\dx
  \geq
   \Big(\frac{n-2}{2}\Big)^2\int_{\Om}\frac{|u|^2}{|x|^2}~\dx,
\ee
with the best possible constant. As in \cite[pg 262]{S}, an integration by parts shows that
\[
 (n-2)\int_{\Om}\frac{|u|^2}{|x|^2}~\dx
  =
   -2\int_{\Om} |u||x|^{-2}x\cdot\nabla|u|~\dx,
\]
and then applying the Cauchy-Schwarz inequality on the right gives (\ref{Hardy:L^2}). The best constant in this easily obtained functional inequality has applications in various branches of analysis. For instance, it is used in \cite[Appendix B]{S} to prove the non-existence of stable minimal cones in $\mathbb{R}^{n+1}$, $2\leq n\leq6$ (see \cite{CHS}, \cite{CC}, \cite{HHS} and also \cite[\S3.4.1]{DHTr} for related non-stability results using the best constant in (\ref{Hardy:L^2})). More widespread is the critical role it plays in the analysis of the heat equation involving the singular potential $1/|x|^2$ (see for example \cite{BrV}, \cite{GP},\cite{Gk1, Gk2} \cite{CM}, \cite{VzZ}, , \cite{FMT}, \cite{FT}, \cite{VnZ} \& \cite{Erv}).

\medskip

It was first shown in \cite{BrV} that (\ref{Hardy:L^2}) can be improved in case $\Om$ has finite Lebesgue measure $\mathcal{L}^n(\Om)$, by adding the term
 \be\label{Hardy-Sobolev:subcritical:L^2}
  C(n,q)\big(\mathcal{L}^n(\Om)\big)^{1/2^\star-1/q}\bigg(\int_{\Om} |u|^q~\dx\bigg)^{2/q},~~~~1\leq q<2^\star:=\frac{2n}{n-2},
 \ee
on its right hand side. Such a non-homogeneous improvement fails in the case of the critical Sobolev exponent $q=2^\star$. This can be seen for example from the improvement given for bounded $\Om$ in \cite[Theorem A]{FT}. There, the term
 \be\label{Hardy-Sobolev:L^2}
  C(n)\bigg(\int_{\Om}|u|^{2^\star}X_1^{1+2^\star/2}(|x|/D_0)~\dx\bigg)^{2/2^\star},~~~~X_1(t):=(1-\log t)^{-1},~t\in(0,1],
 \ee
where $D_0:=\sup_{x\in\Om}|x|$, was added on the right hand side of (\ref{Hardy:L^2}). It is then shown that the exponent $1+2^\star/2$ on $X_1$ cannot be decreased, stating thus the failure of adding the term (\ref{Hardy-Sobolev:subcritical:L^2}) for $q=2^\star$.

It was also established in \cite[Theorem D]{FT} in case of bounded $\Om$, that the homogeneous term
\be\label{Hardy-improved:L^2}
  \frac{1}{4}\int_{\Om}\frac{|u|^2}{|x|^2}X_1^2(|x|/D_0)~\dx,
 \ee
can be added on the right hand side of (\ref{Hardy:L^2}). Here, the constant $1/4$ is the best possible and the exponent $2$ on $X_1$ cannot be decreased. This type of optimal logarithmic homogeneous improvements to Hardy inequalities originated in \cite{BrM}. Compared to (\ref{Hardy-Sobolev:L^2}), one has that (\ref{Hardy-improved:L^2}) is not a weaker remainder term. In particular, one cannot deduce (\ref{Hardy:L^2}) with remainder term (\ref{Hardy-improved:L^2}) from (\ref{Hardy:L^2}) with remainder term (\ref{Hardy-Sobolev:L^2}) by applying H\"{o}lder's inequality (not even with some other positive constant instead of $1/4$). However, the optimality of the exponent 2 on $X_1$ in (\ref{Hardy-improved:L^2}) implies via H\"{o}lder's inequality the optimality of the exponent $1+2^\star/2$ on $X_1$ in (\ref{Hardy-Sobolev:L^2}) (see \cite{PsSp} and \cite{FPs} for similar arguments).

In the question what is an optimal non-homogeneous term one can add on the right hand side of (\ref{Hardy:L^2}) with remainder term (\ref{Hardy-improved:L^2}), the answer is
\[
 C(n)\bigg(\int_{\Om}|u|^{2^\star}X_1^{1+2^\star/2}(|x|/D_0) X_2^{1+2^\star/2}(|x|/D_0)~\dx\bigg)^{2/2^\star},~~~~X_2(t):=X_1(X_1(t)),~t\in(0,1].
\]
Furthermore, it is proved in \cite[Theorem A$^\prime$]{FT} that for any $k\in\N\cup\{0\}$ and all $u\in H^1_0(\Om)$
 \begin{align}\nonumber &
  \int_{\Om}|\nabla u|^2~\dx
   -
    \Big(\frac{n-2}{2}\Big)^2\int_{\Om}\frac{|u|^2}{|x|^2}~\dx
     -
      \frac{1}{4}\int_{\Om}\frac{|u|^2}{|x|^2}\sum_{i=1}^k\prod_{j=1}^iX_j^2(|x|/D_0)~\dx
       \\ \label{Hardy-Sobolev-k:L^2} & \geq
        C(n)\bigg(\int_{\Om}|u|^{2^\star}\prod_{i=1}^{k+1}X_i^{1+2^\star/2}(|x|/D_0)~\dx\bigg)^{2/2^\star},
         ~~~~X_{j+1}(t):=X_1(X_j(t)),~t\in(0,1],
 \end{align}
where the exponent $1+2^\star/2$ on $X_{k+1}$ cannot be decreased. That the left hand side in this inequality is nonnegative and each term appears with best constant $1/4$ and optimal exponent $2$ on $\prod_{j=i}^kX_j$ was first established in \cite[Theorem D]{FT}. For a second proof of (\ref{Hardy-Sobolev-k:L^2}) with the best constant $C(n)$ see \cite{AdFT}.

\medskip

The purpose of this paper is to extend inequality (\ref{Hardy-Sobolev-k:L^2}) to the case of the $k$-improved $p$-Hardy inequality for any $1<p<\infty$, $p\neq n$. More precisely, the $p$-Hardy inequality in a domain $\Om$ of $\rn$, $n\geq2$, containing the origin, asserts that if $p\geq1$, $p\neq n$, then for all $u\in W_0^{1,p}(\Om)$ if $p<n$, or all $u\in W_0^{1,p}(\Om\setminus\{0\})$ if $p>n$, we have
 \be\label{Hardy}
  \int_{\Om}|\nabla u|^p~\dx
   \geq
    \Big|\frac{n-p}{p}\Big|^p\int_{\Om}\frac{|u|^p}{|x|^p}~\dx,
 \ee
with the best possible constant. A proof of the same simplicity as in the case $p=2$ follows by the analogous integration by parts and H\"{o}lder's inequality. It is proved in \cite{BFT2} by a vector field method, that in bounded domains $\Om$ containing the origin there exists $b=b(n,p)\geq1$ such that for any $k\in\mathbb{N}$ the terms
 \be\label{Hardy-k-improved:L^p}
  \frac{p-1}{2p}\Big|\frac{n-p}{p}\Big|^{p-2}\int_{\Om}\frac{|u|^p}{|x|^p}\sum_{i=1}^k\prod_{j=1}^iX_j^2(|x|/D)~\dx,~~~~D=bD_0,
 \ee
can be added on its right hand side. Moreover, for every $k\in\mathbb{N}$, each one of these terms appears with the best possible constant and the exponent $2$ on $\prod_{j=1}^kX_j$ cannot be decreased. In Theorem \ref{BFT:Theorem:improved} we provide an alternative proof of \eqref{Hardy-k-improved:L^p} using a suitable ground state transformation (as was done in the $p=2$ case in \cite[Theorem D]{FT}).

\medskip

Our first result is the extension to all values of $p\in(1,n)$ of \eqref{Hardy-Sobolev-k:L^2}. We denote below by $p^\star$ the critical Sobolev exponent in this case; that is \[p^\star:=np/(n-p).\]

\bigskip \noindent {\bf Theorem A} {\em Let $\Om$ be a bounded domain in $\rn$, $n\geq2$, containing the origin and let $1<p<n$. There exist constants $B=B(n,p)\geq1$ and $C=C(n,p)>0$ such that for any $k\in\N\cup\{0\}$ and all $u\in W_0^{1,p}(\Om)$
 \begin{align}\nonumber &
  \int_{\Om}|\nabla u|^p~\dx
   -
    \Big(\frac{n-p}{p}\Big)^p\int_{\Om}\frac{|u|^p}{|x|^p}~\dx
     -
      \frac{p-1}{2p}\Big(\frac{n-p}{p}\Big)^{p-2}\int_{\Om}\frac{|u|^p}{|x|^p}\sum_{i=1}^k\prod_{j=1}^iX_j^2(|x|/D)~\dx
       \\ \label{Hardy-Sobolev-k:L^p<n} & \geq
        C\bigg(\int_{\Om}|u|^{p^\star}\prod_{i=1}^{k+1}X_i^{1+p^\star/p}(|x|/D)~\dx\bigg)^{p/p^\star},
        ~~~~D=B\sup_{x\in\Om}|x|.
 \end{align}
Moreover, for each $k\in\N\cup\{0\}$, the weight function $\prod_{i=1}^{k+1}X_i^{1+p^\star/p}$ is optimal in the sense that the power $1+p^\star/p$ on $X_{k+1}$ cannot be decreased.}

\bigskip

\noindent Inequality \eqref{Hardy-Sobolev-k:L^p<n} for $p\neq2$ is established here for the first time even with $k=0$. In fact, for $k=0$ the subcritical ($L^q$-weighted remainder term with $q<p^\star$) inequality with an optimal power on the logarithmic weight, that tends to the one appearing in \eqref{Hardy-Sobolev-k:L^p<n} as $q\rightarrow p^\star$, is in \cite[Theorem C (2)]{BFT1}. It is referred there as an open question whether (\ref{Hardy-Sobolev-k:L^p<n}) for $k=0$ was true. Prior results in the direction of obtaining \eqref{Hardy-Sobolev-k:L^p<n} with $k=0$ include \cite[Theorem 6.4]{BFT1}, \cite[Theorem 1.1]{AdChR} and \cite[Theorem 1.1]{AbdCP}.

\medskip

The proof of Theorem A splits in two cases regarding to whether $p<2$ or $p>2$:

\medskip

For $1<p<2$ we push further the basic idea developed in \cite{AdFT}. In particular, applying first a suitable ground state transformation; cf. \cite{BrV}, we get that \eqref{Hardy-Sobolev-k:L^p<n} will follow once a critical case of the Caffarelli-Kohn-Nirenberg inequality is true. To establish this inequality we compare it's best constant with the best Sobolev constant by applying an Emden-Fowler transform; cf. \cite{CW}. A notable difficulty in the proof is to establish Theorem \ref{BFT:Theorem:improved} $(ii)$ which is the extension of Proposition 3.4 of \cite{BFT1} to all $k\in\N$. The reason is that for $k\in\N$, the naturally choice for a function to be used in the ground state transformation fails to be a supersolution of the corresponding Euler-Langrange equation (see also Remark \ref{rmrkGSHp<2}). To solve this problem we invent a new supersolution (see \eqref{fak}).

\medskip

For $p>2$, by applying the natural ground state transform, two lower bounds for the left hand side of \eqref{Hardy-Sobolev-k:L^p<n} may be produced (for $k=0$ see \cite[Proposition 3.3 eq. (3.6) \& (3.7)]{BFT1}). Introducing a third lower bound simply by adding these two lower bounds, we show that decomposition in spherical harmonics (see \cite{VzZ}) works for $p\neq 2$ as well (see \cite{FT} for $p=2$)! Still, each lower bound has to be treated separately before they are added. By a further change of variables one of these bounds is reduced to a one dimensional integral which yields the correct exponent on the logarithmic weight. Using a new argument, where the Poincar\'{e} inequality on the sphere comes into play, we manage to get a cooperative estimate for the second lower bound.\\

Passing to the case $p>n$, we address here the question of what is an optimal nonhomogeneous term one can add on the right hand side of (\ref{Hardy}) with remainder term (\ref{Hardy-k-improved:L^p}). It is known that the Dirichlet integral in $\rn$ exceeds a constant multiple of the $C^{0,1-n/p}$-seminorm. More precisely,
there exists a positive constant $C=C(n,p)$ such that for all $u\in W^{1,p}(\rn)$
 \be\nonumber\label{Morrey:Sobolev:p>n}
   \bigg(\int_{\rn}|\nabla u|^pdx\bigg)^{1/p}
    \geq C
     \sup_{\substack{x,y\in\rn \\ x\neq y}}
      \frac{|u(x)-u(y)|}{|x-y|^{1-n/p}},
 \ee
and the modulus of continuity $1-n/p$ is optimal. This is Morrey's inequality. In Hardy's inequality \eqref{Hardy}, an optimally weighted $C^{0,1-n/p}$-seminorm was added in \cite{Ps} in case of a bounded $\Om$. The precise statement asserts the existence of constants $C=C(n,p)>0$ and $B=B(n,p)\geq1$ such that for all $u\in W_0^{1,p}(\Om\setminus\{0\})$
 \be\nonumber\label{Hardy:Morrey:Sobolev:p>n}
  \bigg(\int_{\Om}|\nabla u|^pdx -\Big(\frac{p-n}{p}\Big)^p \int_{\Om}\frac{|u|^p}{|x|^p}dx\bigg)^{1/p}
   \geq C
    \sup_{\substack{x,y\in\Om \\ x\neq y}}
     \frac{|u(x)-u(y)|}{|x-y|^{1-n/p}}X_1^{1/p}\Big(\frac{|x-y|}{D}\Big),~~~~D=B\diam(\Om).
 \ee
The correction $X_1^{1/p}$ on the modulus of continuity was shown to be optimal in the sense that the power $1/p$ on $X_1$ cannot be decreased. The following inequality is reduced to the above one when $k=0$, and gives the complete picture for the series improvement of Hardy's inequality

\bigskip \noindent {\bf Theorem B} {\em Let $\Om$ be a bounded domain in $\rn$, $n\geq2$, containing the origin and let $p>n$. There exist constants $B=B(n,p)\geq1$ and $C=C(n,p)>0$ such that for any $k\in\N\cup\{0\}$ and all $u\in W_0^{1,p}(\Om\setminus\{0\})$
 \begin{align}\nonumber &
  \bigg(\int_{\Om}|\nabla u|^p~\dx
   -
    \Big(\frac{p-n}{p}\Big)^p\int_{\Om}\frac{|u|^p}{|x|^p}~\dx
     -
      \frac{p-1}{2p}\Big(\frac{p-n}{p}\Big)^{p-2}\int_{\Om}\frac{|u|^p}{|x|^p}\sum_{i=1}^k\prod_{j=1}^iX_j^2(|x|/D)~\dx\bigg)^{1/p}
       \\ \label{k:improvedHardy:Morrey} & \geq
        C\sup_{\substack{x,y\in\Om \\ x\neq y}}
        \frac{|u(x)-u(y)|}{|x-y|^{1-n/p}}\prod_{i=1}^{k+1}X_i^{1/p}(|x-y|/D),~~~~D=B\diam(\Om).
 \end{align}
Moreover, for each $k\in\N\cup\{0\}$, the weight function $\prod_{i=1}^{k+1}X_i^{1/p}$ is optimal in the sense that the power $1/p$ on $X_{k+1}$ cannot be decreased.}

\bigskip

The paper is organised as follows: In \S2, after setting the notation and a couple of technical calculus facts, we use the ground state transform to produce lower estimates for the series expansion of Hardy's inequality. These are used in \S3, \S4 and \S6 to prove Theorem A for $1<p<2$, Theorem A for $p>2$ and Theorem B for $p>n$, respectively. These are also used in \S5 to prove a suitable local estimate on balls that is crucial for the proof of Theorem B.
\section{Preparative results}
In this paper we assume
\begin{itemize}
\item[$\bullet$] $1<p\neq n$, where $n\in\N\setminus\{1\}$,
\item[$\bullet$] $\Om$ is a bounded domain in $\rn$ containing the origin,
\item[$\bullet$] $D_0:=\sup_{x\in\Om}|x|$.
\end{itemize}
Furhermore, $\mathcal{L}^n$ stands for the Lebesgue measure in $\rn$ and $\mathcal{H}^{n-1}$ for the $n-1$ Hausdorff measure in $\rn$. $B_r(x)$ is the open ball in $\rn$ having radius $r>0$ and centre at $x\in\rn$; $\prt B_r(x)$ is its boundary. When the centre is of no importance we simply write $B_r$. When the center is the origin and $r=1$ we write $\us$ instead of $\prt B_1(0)$. Also, $\omn:=\mathcal{L}^n(B_1)$ and so $\mathcal{H}^{n-1}(\prt B_1)=n\omn$. Throughout, an expression of the form $b(n,p,...),~B(n,p,...)$, $c(n,p,...)$ or $C(n,p,...)$ stands for a positive constant that may change value from line to line but always depending only on its arguments $n,p...$. The particular constant depending only on $p$ that appears in \eqref{vecineq} or \eqref{vecineq2}, is denoted by $c_p$. All functions having compact support are extended by zero outside it.
\subsection{Some calculus facts}
\begin{definition}For any $t\in(0,1]$ we define the function $X_1(t):=(1-\log t)^{-1}$ and then
 \[
  X_k(t):=X_1(X_{k-1}(t)),~k=2,3,...,
   ~~~~Y_k(t):=\prod_{i=1}^kX_i(t),~k\in\N,
    ~~~~Z_k(t):=\sum_{i=1}^kY_i(t),~k\in\N.
 \]
\end{definition}

The following computational lemma gives a formula for the derivative of $X_k,~Y_k$ and $Z_k$.

\begin{lemma}\label{lemma:differential:rules} For any $k\in\N$ and $t\in(0,1]$ there holds
 \[
  \frac{d}{dt}\big(X_k(t)\big)=\frac{1}{t}Y_k(t)X_k(t),
   ~~~~\frac{d}{dt}\big(Y_k(t)\big)=\frac{1}{t}Y_k(t)Z_k(t),
    ~~~~\frac{d}{dt}\big(Z_k(t)\big)=\frac{1}{2t}\Big(Z_k^2+\sum_{i=1}^kY_i^2(t)\Big).
 \]
\end{lemma}

\noindent\textbf{Proof.} The first one follows easily by induction. The proof of the second one is
 \[
  \frac{d}{dt}\big(Y_k(t)\big)
   =\sum_{j=1}^k\frac{d}{dt}\big(X_j(t)\big)\prod_{\substack{i=1 \\ i\neq j}}^kX_i(t)
    =\frac{1}{t}\sum_{j=1}^kY_j(t)X_j(t)\prod_{\substack{i=1 \\ i\neq j}}^kX_i(t)
     =\frac{1}{t}Y_k(t)Z_k(t),
 \]
where the first one is used in the middle equality. For the third one, notice that one has the elementary identity
 \[
  \sum_{i=1}^kY_iZ_i
   =
    \frac{1}{2}\Big(Z_k^2+\sum_{i=1}^kY_i^2\Big),
 \]
for which we include its proof for clarity
 \[
   2\sum_{i=1}^kY_iZ_i
   =2\sum_{\substack{i,j=1 \\ j\leq i}}^kY_iY_j
    =2\sum_{\substack{i,j=1 \\ j<i}}^kY_iY_j+2\sum_{i=1}^kY^2_i
     =\Big(\sum_{i=1}^kY_i\Big)^2+\sum_{i=1}^kY^2_i
      =Z_k^2+\sum_{i=1}^kY_i^2.
 \]
Now we may easily conclude
 \[
  \frac{d}{dt}\big(Z_k(t)\big)
    =\sum_{i=1}^k\frac{d}{dt}\big(Y_i(t)\big)
     =\frac{1}{t}\sum_{i=1}^kY_i(t)Z_i(t)
      =\frac{1}{2t}\Big(Z_k^2(t)+\sum_{i=1}^kY_i^2(t)\Big),
 \]
where the second one is used in the middle equality. \qed

\begin{remark}\label{remark:Zinfinite}
The infinite series $Z_\infty(t):=\sum_{k=1}^{\infty}Y_k(t)$, $t\in(0,1]$, converges if and only if $t\in(0,1)$. A proof of this fact can be extracted from \cite[\S 6]{FT} (see \cite[Appendix]{D} for the details).
\end{remark}

\noindent A technical lemma follows

\begin{lemma}\label{Lemma:X}
Let $\al,\beta,R>0$. For all $r\in(0,R]$, all $c>1/\alpha$ and any $D\geq\eta R$, where $\eta=\eta(\al,\beta,c)>1$, there holds
\be\nonumber
 \int_0^rt^{\al-1}Y_k^{-\beta}(t/D)\mathrm{d}t
  \leq
   cr^\al Y_k^{-\beta}(r/D).
\ee
\end{lemma}

\noindent\textbf{Proof.} For $c>0$ and $D\geq R$ set
\be\nonumber
 f(r)
  :=
   \int_0^rt^{\al-1}Y_k^{-\beta}(t/D)\mathrm{d}t
    -
     cr^\al Y_k^{-\beta}(r/D),~~~~r\in(0,R].
\ee
It suffices to show for suitable values of the parameters $c$ and $D$, that $f(r)\leq0$ for all $r\in(0,R)$. Since $f(0+)=0$, it is enough to choose $c$ and $D$ so that $f$ is decreasing in $(0,R)$. To this end, with the aid of Lemma \ref{lemma:differential:rules} we compute
 \be\nonumber
f'(r)
 =
  cr^{\al-1}Y_k^{-\beta}(r/D)\big[1/c-\al+\beta Z_k(r/D)\big],
   ~~~~r\in(0,R].
 \ee
By Remark \ref{remark:Zinfinite} the series $Z_\infty(R/D)$ is convergent if $R<D$, hence for $c>1/\alpha$ we can find large enough $\eta>1$ that depends only on $\al,\beta$ and $c$, such that for all $D\geq\eta R$ to have $f'(r)\leq0$ for all $r\in(0,R)$. \qed
\subsection{Improvements via the ground state transform}
\begin{definition} Given $u\in\test(\Om\setminus\{0\})$ and $D\geq D_0$ we set
\[
 I_0[u;D]\equiv I_0[u]
  :=
   \int_{\Om}|\nabla u|^p~\dx
    -
     \Big|\frac{n-p}{p}\Big|^p\int_{\Om}\frac{|u|^p}{|x|^p}\dx,
\]
\begin{align}\nonumber
 I_k[u;D]
  & :=
   I_{k-1}[u;D]
    -
     \frac{p-1}{2p}\Big|\frac{n-p}{p}\Big|^{p-2}\int_{\Om}\frac{|u|^p}{|x|^p}Y_k^2(|x|/D)\dx
      \\ \nonumber & =
       I_0[u]-\frac{p-1}{2p}\Big|\frac{n-p}{p}\Big|^{p-2}\int_{\Om}\frac{|u|^p}{|x|^p}\sum_{i=1}^kY_i^2(|x|/D)\dx,
        ~~~~k\in\N.
\end{align}
\end{definition}
In \cite{BFT2} the following successive homogeneous improvements to Hardy's inequality were obtained

\begin{theorem}[\cite{BFT2}-Theorem A]\label{BFT:Theorem}
There exists a constant $b=b(n,p)\geq1$ such that for any $k\in\N$
\be\label{improved_Hardy:BFT}
 I_k[u;D]\geq0~~~~\mbox{for all }u\in\test(\Om\setminus\{0\}),
\ee
where $D=bD_0$. Moreover, for each $k\in\N$:
\begin{itemize}
\item[$(i)$] the weight function $Y_k^2$ is optimal, in the sense that the power 2 cannot be decreased, and
\item[$(ii)$] the constant appearing on the term with this weight function is sharp.
\end{itemize}
\end{theorem}

In \cite{BFT1} the authors obtained various auxiliary improvements for Hardy's inequality (\ref{Hardy}). In particular, given $u\in\test(\Om\setminus\{0\})$, the ground state transformation
\be\label{GSH}
 u(x)=|x|^{1-n/p}v(x),
\ee
plus elementary vectorial inequalities lead to the following lower bounds on $I_0[u]$ in terms of the function $v$ (see \cite[Lemma 3.3 \& Proposition 3.4]{BFT1})
\[
 I_0[u]\geq c(p)\int_{\Om}|x|^{p-n}|\nabla v|^p~\dx,
\]
\[
 I_0[u]\geq c(p)\int_{\Om}|x|^{2-n}|v|^{p-2}|\nabla v|^2\dx,
\]
both in case $p\geq2$, and
\be\label{GSHp<2}
 I_0[u]\geq c(n,p)\int_{\Om}|x|^{p-n}|\nabla v|^pX_1^{2-p}(|x|/D)\dx,~~~~D\geq D_0,
\ee
in case $p<2$. Our aim here is to extend these estimates to arbitrary $k\in\N$. More precisely, we have the following theorem which readily implies (\ref{improved_Hardy:BFT})

\begin{theorem}\label{BFT:Theorem:improved}
For $a\geq0$, $D\geq D_0$ and $k\in\N\cup\{0\}$ set
\be\label{fak}
 f_{a,k,D}(x):=\sgn(n-p)|x|^{1-n/p}Y_k^{-1/p}(|x|/D)\big(1-aX_1(|x|/D)\big),~~~~x\in\Om\setminus\{0\}.
\ee
For simplicity we write $f_0$ in place of $f_{0,0,D}$ and $f_{k,D}$ instead of $f_{0,k,D}$. Then,
\begin{itemize}
\item[$(i)$] for $p\geq2$, there exists a constant $b'=b'(n,p)\geq1$ such that for any $k\in\N\cup\{0\}$, all $u\in\test(\Om\setminus\{0\})$ and any $D\geq b'D_0$
\be\label{BFTprime:p>2}
 I_k[u;D]
  \geq
   c_p\int_{\Om}|x|^{p-n}|\nabla v|^pY_k^{-1}(|x|/D)\dx,
\ee
\be\label{BFTprime:p>2:Sobolev}
 I_k[u;D]
  \geq
   c_p
    \int_{\Om}|x|^{2-n}|v|^{p-2}|\nabla v|^2\Big|\frac{p-n}{p}-\frac{1}{p}Z_k(|x|/D)\Big|^{p-2}Y_k^{-1}(|x|/D)\dx,
\ee
where $v$ is defined through the ground state transformation $u=f_{k,D}v$.
\item[$(ii)$] If $p<2$, there exist constants $a=a(n,p)>0$ and $b''=b''(n,p)\geq1$ such that for any $k\in\N\cup\{0\}$, all $u\in\test(\Om\setminus\{0\})$ and any $D\geq b''D_0$
\be\label{BFTprime:p<2}
 I_k[u;D]
  \geq
   c(n,p)\int_{\Om}|x|^{p-n}|\nabla v|^pY_{k+1}^{2-p}(|x|/D)Y_k^{-1}(|x|/D)\dx,
\ee
where $v$ is defined through the ground state transformation $u=f_{a,k,D}v$.
\end{itemize}
\end{theorem}

\begin{remark}\label{rmrk_convergence}
Clearly,
\[p<n
 \Longrightarrow \Big|\frac{p-n}{p}-\frac{1}{p}Z_k(|x|/D)\Big|=\frac{n-p}{p}+\frac{1}{p}Z_k(|x|/D)\geq\frac{n-p}{p},
\]
for all $x\in\Om$ and all $k\in\N$. A similar estimate is true for $p>n$. In particular, by Remark \ref{remark:Zinfinite} the series $Z_\infty(|x|/D_0)$ is convergent if $|x|<D_0$ and thus for suitable $b'''>1$, depending only on $n,p$, we may choose $D\geq b'''D_0$ so that
\[
 \frac{p-n}{p}-\frac{1}{p}Z_k(|x|/D)\geq C(n,p)>0,
\]
for all $x\in\Om$ and all $k\in\N$. Consequently, from \eqref{BFTprime:p>2:Sobolev} we get
\end{remark}

\begin{corollary}
For $p\geq2$, there exists a constant $b'''=b'''(n,p)\geq1$ such that for all $u\in\test(\Om\setminus\{0\})$, any $D\geq b'''D_0$ and any $k\in\N\cup\{0\}$
\be\label{BFTdoubleprime:p>2:Sobolev}
 I_k[u;D]
  \geq
   C(n,p)
    \int_{\Om}|x|^{2-n}|v|^{p-2}|\nabla v|^2Y_k^{-1}(|x|/D)\dx.
\ee
\end{corollary}

To prove Theorem \ref{BFT:Theorem:improved} we need the following key lemma

\begin{lemma}\label{super-sub-sol}\begin{itemize}
\item[$(i)$] For $p\geq2$, there exists a constant $b'=b'(n,p)\geq1$ such that for any $D\geq b'D_0$ and any $k\in\N$, the function $f_{k,D}$ (defined by \eqref{fak} with $a=0$) is a supersolution of the following $p$-Laplace equation with $k+1$ singular potential terms
 \be\label{plaplace:infseries}
  -\Delta_pw
   -\Big(\Big|\frac{p-n}{p}\Big|^p+\frac{p-1}{2p}\Big|\frac{p-n}{p}\Big|^{p-2}\sum_{i=1}^kY_i^2(|x|/D)\Big)\frac{|w|^{p-2}w}{|x|^p}
    =0, ~~~~\mbox{in }\Om\setminus\{0\}.
 \ee
\item[$(ii)$] For $p<2$, there exist constants $a=a(n,p)>0$ and $b''=b''(n,p)\geq1$ such that for any $D\geq b''D_0$ and any $k\in\N$, the function $f_{a,k,D}$ (defined in \eqref{fak}) is a supersolution of \eqref{plaplace:infseries}.
\end{itemize}

\end{lemma}

\noindent\textbf{Proof.} Let $a\geq0$ and $0<\e<1$. In view of Remark \ref{remark:Zinfinite} we choose $\delta=\delta(a,p)\geq1$, such that with $D:=D_0\delta$ to have
\be\label{fragma}
1-aX_1(|x|/D)\geq2-p\quad\text{and}\quad X_1(|x|/D)\leq\sum_{i=1}^\infty Y_i(|x|/D)\leq pX_1(|x|/D),~~\forall~x \in \Om.
\ee
We further set
\be\label{def:Ak}
 A_{a,k}(|x|/D)
  :=
   \frac{p-n}{p}-\frac{1}{p}Z_k(|x|/D) -a\frac{X_1^{2}(|x|/D)}{1-aX_1(|x|/D)}.
\ee
Using Lemma \ref{lemma:differential:rules} we compute (from now on in this proof we write $f_k$, $A_k$, $X_k$, $Y_k$, $Z_k$ instead of $f_{a,k,D}(x)$, $A_{a,k}(|x|/D)$, $X_k(|x|/D)$, $Y_k(|x|/D)$, $Z_k(|x|/D)$)
\[
 \nabla f_k=\frac{f_k}{|x|}A_k\frac{x}{|x|}
  ~~~~\mbox{so that}~~~~
   -\Delta_pf_k=-\Div\Big\{\frac{|f_k|^{p-2}f_k}{|x|^{p-1}}|A_k|^{p-2}A_k\frac{x}{|x|}\Big\}.
\]
Direct computations reveal the next identities which are valid for any $x\in\Om\setminus\{0\}$
\[
 -\Div\Big\{\frac{|f_k|^{p-2}f_k}{|x|^{p-1}}\frac{x}{|x|}\Big\}
  =
   \Big(p-n-(p-1)A_k\Big)\frac{|f_k|^{p-2}f_k}{|x|^p},
\]
\begin{align}\nonumber
 -\nabla\Big(|A_k|^{p-2}A_k\Big)
   & =
   -(p-1)|A_k|^{p-2}\nabla A_k
    \\ \nonumber & =
     \frac{p-1}{|x|}|A_k|^{p-2}\bigg(\frac{1}{2p}\Big(Z^2_k+\sum_{i=1}^kY^2_i\Big) +\underbrace{2a\frac{X_1^{3}}{1-aX_1}+a^2\frac{X_1^{4}}{(1-aX_1)^2}}_{:=F(X_1)}\bigg)\frac{x}{|x|},
\end{align}
where in the last one we used Lemma \ref{lemma:differential:rules} in order to compute $\nabla A_k$. We conclude
\begin{align}\nonumber
      -\Delta_pf_k & =
  |A_k|^{p-2}
  \bigg((p-n)A_k-(p-1)A_k^2+\frac{p-1}{2p}\Big(Z^2_k+\sum_{i=1}^kY^2_i\Big)+(p-1)F(X_1)\bigg)
  \frac{|f_k|^{p-2}f_k}{|x|^p}
      \\ \nonumber & =
  |A_k|^{p-2}
  \bigg(\Big(\frac{p-n}{p}\Big)^2+\frac{(p-n)(p-2)}{p^2}Z_k+\frac{(p-1)(p-2)}{2p^2}Z^2_k+\frac{p-1}{2p}\sum_{i=1}^kY^2_i\bigg)
  \frac{|f_k|^{p-2}f_k}{|x|^p}\\ \nonumber
      & ~~~~ +
  a|A_k|^{p-2}
  \bigg(\frac{(p-n)(p-2)}{p}+2(p-1)\Big(X_1-\frac{1}{p}Z_k\Big)\bigg)
  \frac{X_1^{2}}{1-aX_1}\frac{|f_k|^{p-2}f_k}{|x|^p}.
\end{align}
It turns out that given $1<p<n$ it is enough to establish the following inequality for some nonnegative constant $a=a(n,p)$ and for any $x\in\Om\setminus\{0\}$
\begin{align} \nonumber &
  |A_k|^{p-2}
  \bigg(\Big(\frac{p-n}{p}\Big)^2+\frac{(p-n)(p-2)}{p^2}Z_k+\frac{(p-1)(p-2)}{2p^2}Z^2_k+\frac{p-1}{2p}\sum_{i=1}^kY^2_i\bigg)
    \\ \nonumber
      & ~~~~ +
  a|A_k|^{p-2}
  \bigg(\frac{(p-n)(p-2)}{p}+2(p-1)\Big(X_1-\frac{1}{p}Z_k\Big)\bigg)
  \frac{X_1^{2}}{1-aX_1}
   \\ \label{suffices} & ~~~~~~~~ \geq
    \Big|\frac{p-n}{p}\Big|^p+\frac{p-1}{2p}\Big|\frac{p-n}{p}\Big|^{p-2}\sum_{i=1}^kY_i^2,
\end{align}
and the reverse inequality if $p>n$ (note that $\sgn f_k=\sgn(n-p)$). In the case $p=2$ we take $a=0$ and this inequality is immediately true since we have equality. In what follows we assume $p\neq2$.

To write (\ref{suffices}) in a more accessible form let us set
\[
 h:=\frac{p-n}{p},~~~~t:=\frac{Z_k}{p},~~~~s:=a\frac{X_1^{2}}{1-aX_1}
 ~~~~\mbox{ and }~~~~ \lambda:=\frac{p-1}{2p}\sum_{i=1}^kY^2_i,
\]
so that all we need to prove is that for all sufficiently small $t$ and some nonnegative constant $a$ depending possibly only on $n,p$, there holds
\[
 \Big|1-\frac{t+s}{h}\Big|^{p-2}
  \Big(h^2+\lambda+h(p-2)(t+s)+\frac{(p-1)(p-2)}{2}t^2+2(p-1)(X_1-t)s\Big)
   \geq
     h^2+\lambda,
\]
if $1<p<n$, and the reverse inequality if $p>n$. By a further rearrangement of terms, this is the same as
\be\label{sufficesLemma2.10}
 (h^2+\lambda)\bigg(1-\Big|1-\frac{t+s}{h}\Big|^{2-p}\bigg)+h(p-2)(t+s)+\frac{(p-1)(p-2)}{2}t^2+2(p-1)(X_1-t)s
  \geq
   0,
\ee
if $1<p<n$, and the reverse inequality if $p>n$. The Taylor expansion of $g(x):=|1-x|^{2-p}$ at $x=0$ is
\[
 g(x)=1+(p-2)x+\frac{(p-1)(p-2)}{2}x^2+\frac{p(p-1)(p-2)}{6}x^3+O(x^4),
\]
and after an easy computation we get that \eqref{sufficesLemma2.10} is equivalent to
\begin{align}\nonumber &
2(p-1)(X_1-t)s+\frac{(p-1)(p-2)}{2}\big(t^2-(t+s)^2\big)
 \\ \label{sufficesLemma2.10prime} &
  -\frac{p-2}{h}\Big(\frac{p(p-1)}{6}(t+s)^3+\lambda(t+s)\Big)+O\big(\lambda(t+s)^2\big)
   \geq 0,
\end{align}
if $1<p<n$, and the reverse inequality if $p>n$.

\medskip

We distinguish two cases:

\medskip

$\bullet$ $2<p\neq n$. In this case we take $a=0$ and hence $s=0$. The trivial fact that $\sum_{i=1}^kY^2_i\leq Z_k^2$, $k\in\N$, is translated to $2\lambda\leq p(p-1)t^2$. Therefore if $p<n$ we have $h<0$ and so there holds
\begin{align}\nonumber &
 -\frac{p-2}{h}\Big(\frac{p(p-1)}{6}t^3+\lambda t\Big)+O\big(\lambda t^2\big)
   \geq 0,
\end{align}
if $2<p<n$, while if $p>n$ then $h>0$ and so the reverse inequality holds true.

\medskip

$\bullet$ $1<p<2$. By \eqref{fragma} we obtain $X_1-t=X_1-Z_k/p>0$ for all $k\in\N$ and we can choose for example $a=p$ and hence $s>0$ such that the first line in \eqref{sufficesLemma2.10prime} is positive. \qed\\

We also need the following elementary Hardy inequality

\begin{lemma}\label{anilog}
For all $w\in W_0^{1,p}(\Om)$ and all $D\geq D_0$, there holds
\[
 \int_{\Om}\frac{|w|^p}{|x|^n}Y_k(|x|/D)X^2_{k+1}(|x|/D)\dx
  \leq p^p\int_{\Om}|x|^{p-n}|\nabla w|^pY_k^{1-p}(|x|/D)X^{2-p}_{k+1}(|x|/D)\dx.
\]
\end{lemma}

\noindent\textbf{Proof.} A direct computation using Lemma \ref{lemma:differential:rules} shows that
 \[
  \Div\big\{|x|^{-n}X_{k+1}(|x|/D)x\big\}
   =
    |x|^{-n}Y_k(|x|/D)X^2_{k+1}(|x|/D),
     ~~~~x\in\Om\setminus\{0\}.
 \]
Hence, integrating by parts,
\begin{align*}
 &\int_{\Om}\frac{|w|^p}{|x|^n}Y_k(|x|/D)X^2_{k+1}(|x|/D)\dx
  = -p\int_{\Om} X_{k+1}(|x|/D)|w|^{p-1}\nabla|w|\cdot\frac{x}{|x|^n}\dx
   \\ & \leq p\int_{\Om} X_{k+1}(|x|/D)|w|^{p-1}|\nabla w||x|^{1-n}\dx\\
    & = p\int_{\Om}\Big\{\frac{|w|^{p-1}}{|x|^{n(1-1/p)}}Y_k^{1-1/p}(|x|/D)X_{k+1}^{2(1-1/p)}(|x|/D)\Big\}
                      \Big\{\frac{|\nabla w|}{|x|^{n/p-1}}Y_k^{1/p-1}(|x|/D)X_{k+1}^{2/p-1}(|x|/D)\Big\}\dx\\
     &\leq p\bigg(\int_{\Om}\frac{|w|^p}{|x|^n}Y_k(|x|/D)X^2_{k+1}(|x|/D)\dx\bigg)^{1-1/p}
                \bigg(\frac{|\nabla w|^p}{|x|^{n-p}}Y_k^{1-p}(|x|/D)X_{k+1}^{2-p}(|x|/D)\dx\bigg)^{1/p}.
\end{align*}
The result follows by rearranging terms and taking the $p$-th power.\qed\\

\noindent\textbf{Proof of Theorem \ref{BFT:Theorem:improved} for $p\geq2$.} Setting $u(x)=f_{k,D}(x)v(x)$ we get
\begin{align}\nonumber\hspace{-1em}
 \int_{\Om}|\nabla u|^p~\dx
  & =
   \int_{\Om}|v\nabla f_{k,D}+f_{k,D}\nabla v|^p~\dx
    \\ \nonumber & \geq
     \int_{\Om}|v|^p|\nabla f_{k,D}|^p~\dx
      +
       c_p\int_{\Om}|f_{k,D}|^p|\nabla v|^p~\dx
        +
         \int_{\Om} f_{k,D}|\nabla f_{k,D}|^{p-2}\nabla f_{k,D}\cdot\nabla|v|^p~\dx,
\end{align}
where we have used the following inequality, valid for all $\alpha,\beta\in\rn$, $n\geq1$ and $p\geq2$
\be\label{vecineq}
 |\alpha+\beta|^p\geq|\alpha|^p+c_p|\beta|^p+p|\alpha|^{p-2}\alpha\cdot\beta.
\ee
Noting that
\be\label{divid}
 \Div\{f_{k,D}|\nabla f_{k,D}|^{p-2}\nabla f_{k,D}\}=|\nabla f_{k,D}|^p+f_{k,D}\Delta_pf_{k,D},
\ee
we perform an integration by parts in the last term to arrive at
\begin{align}\nonumber
 \int_{\Om}|\nabla u|^p~\dx
  & \geq
   c_p\int_{\Om}|f_{k,D}|^p|\nabla v|^p~\dx
    -
     \int_{\Om}|v|^pf_{k,D}\Delta_pf_{k,D}~\dx
      \\ \nonumber & =
       c_p\int_{\Om}|f_{k,D}|^p|\nabla v|^p~\dx
        -
         \int_{\Om}|u|^pf_{k,D}^{-1}|f_{k,D}|^{2-p}\Delta_pf_{k,D}~\dx.
\end{align}
Inequality (\ref{BFTprime:p>2}) follows now from Lemma \ref{super-sub-sol}-$(i)$. If instead of (\ref{vecineq}) we use
\be\label{vecineq2}
 |\alpha+\beta|^p\geq|\alpha|^p+c_p|\alpha|^{p-2}|\beta|^2+p|\alpha|^{p-2}\alpha\cdot\beta,
\ee
valid for all $\alpha,\beta\in\rn$, $n\geq1$ and $p\geq2$, we similarly obtain (\ref{BFTprime:p>2:Sobolev}) from Lemma \ref{super-sub-sol}-$(i)$.\qed\\

\noindent\textbf{Proof of Theorem \ref{BFT:Theorem:improved} for $1<p<2$.} By the fact that (see \cite{L})
\[
 |\alpha+\beta|^p-|\alpha|^p \geq p|\alpha|^{p-2}\alpha\cdot\beta + \frac{3p(p-1)}{16}\frac{|\beta|^2}{(|\alpha|+|\beta|)^{2-p}},
\]
for all $\alpha,\beta\in\rn$, $n\geq1$ and $p\in(1,2)$, we get setting $u(x)=f_{a,k,D}(x)v(x)$,
\begin{align}\nonumber
 \int_{\Om}|\nabla u|^p~\dx
  & =
   \int_{\Om}|v\nabla f_{a,k,D} +f_{a,k,D}\nabla v|^p~\dx
    \\ \nonumber & \geq
     \int_{\Om}|v|^p|\nabla f_{a,k,D}|^p~\dx
      + \int_{\Om} f_{a,k,D}|\nabla f_{a,,k,D}|^{p-2}\nabla f_{a,,k,D}\cdot\nabla|v|^p~\dx
       \\ \nonumber & ~~ + c(p)\int_{\Om}\frac{f_{a,k,D}^2|\nabla v|^2}{\big(|v||\nabla f_{a,k,D}|+|f_{a,k,D}||\nabla v|\big)^{2-p}}\dx.
\end{align}
By the same integration by parts and \eqref{divid}, but this time using Lemma \ref{super-sub-sol}-$(ii)$, we get for any $D\geq b'D_0$
\be\label{ani1}
 I_k[u;D]
  \geq
   c(p)\int_{\Om}\frac{f_{a,k,D}^2|\nabla v|^2}{\big(|v||\nabla f_{a,k,D}|+|f_{a,k,D}||\nabla v|\big)^{2-p}}~\dx=:c(p)M_2.
\ee
Next we define $M_4$ to have the same integrand as in $M_2$ but with the measure $\rho^{-p}~\dx$ in place of $\dx$, where $\rho(x):=1-aX_1(|x|/D)$, $x\in\Om$. Also, we set
\[
 M_1:=\int_{\Om}f^p_{k,D}|\nabla v|^pY_{k+1}^{2-p}~\dx,
  ~~~~~~~~
 M_3:=\int_{\Om}\frac{f^p_{k,D}}{|x|^p}|v|^pY_{k+1}^{2}\dx.
\]
To get \eqref{BFTprime:p<2} from \eqref{ani1}, it suffices to show $M_1\leq C(n,p)M_2$. Noting that $f_{k,D}=f_{a,k,D}/\rho$, we use H\"{o}lder's and Minkowski's inequalities as follows
\begin{align}\nonumber
 M_1
  & =
   \int_{\Om}\frac{f_{a,k,D}^p|\nabla v|^p}{\big(|v||\nabla f_{a,k,D}|+|f_{a,k,D}||\nabla v|\big)^{(2-p)p/2}}
                    \big(|v||\nabla f_{a,k,D}|+|f_{a,k,D}||\nabla v|\big)^{(2-p)p/2}Y_{k+1}^{2-p}~\frac{\dx}{\rho^p}
    \\ \nonumber & \leq
      M_4^{p/2}\Big(\int_{\Om}\big(|v||\nabla f_{a,k,D}|+|f_{a,k,D}||\nabla v|\big)^pY_{k+1}^{2}\frac{\dx}{\rho^p}\Big)^{1-p/2}
       \\ \label{ani2} & \leq
         M_4^{p/2}\bigg(\Big(\int_{\Om} |v|^p|\nabla f_{a,k,D}|^pY_{k+1}^{2}\frac{\dx}{\rho^p}\Big)^{1-p/2}+M_1^{1-p/2}\bigg),
\end{align}
where we also have used $(\alpha+\beta)^q\leq\alpha^q+\beta^q$ for all $\alpha,\beta\geq0$, $q=p(1-p/2)\in(0,1]$ and the simple fact that $Y^2_{k+1}(t)\leq Y^{2-p}_{k+1}(t)$ for all $t\in(0,1]$. From Remark \ref{remark:Zinfinite} we know that for sufficiently large $D\geq BD_0$, $B=B(n,p)\geq1$, we have $Z_k(|x|/D)\leq C(n,p)$ for all $x\in\Om$. Hence, taking into account \eqref{fragma} we get $|A_{a,k}(|x|/D)|\leq C(n,p)$ for all $x\in\Om$, where $A_{a,k}$ is given by \eqref{def:Ak}. Therefore
\[
 |\nabla f_{a,k,D}|
  =
   \Big|\frac{f_{a,k,D}}{|x|}A_{a,k}\frac{x}{|x|}\Big|
    \leq
     C(n,p)\frac{|f_{a,k,D}|}{|x|},
\]
so that
\[
 \int_{\Om} |v|^p|\nabla f_{a,k,D}|^pY_{k+1}^{2}\frac{\dx}{\rho^p}
  \leq C(n,p)M_3.
\]
But notice that Lemma \ref{anilog} asserts $M_3\leq p^p M_1$. Plugging these into \eqref{ani2} we arrive at
\[
 M_1\leq C(n,p)M_4.
\]
Finally, $M_4\leq(2-p)^{-p}M_2$ because of \eqref{fragma}, and the proof is complete. \qed

\begin{remark}\label{remark_equality}
In the case $p=2$ the estimate \eqref{BFTprime:p>2} is valid with equality for any $D\geq D_0$ and $B(n,p)=1$.
\end{remark}
\section{Proof of Theorem A when $1<p\leq2$}\label{section:p<2}
We start with a series of reductions. First, since $0\in\Om$, if $u\in\test(\Om)$ then $u\in\test(B_{D_0}(0))$. Hence it is enough to establish \eqref{Hardy-Sobolev-k:L^p<n} for $\Om=B_{D_0}(0)$. Furthermore, \eqref{Hardy-Sobolev-k:L^p<n} being scaling invariant, it is enough to prove it for $\Om=B_1(0)$ only. Finally,  given $u\in\test(B_1(0))\setminus\{0\}$, the transform $u=f_{a,k,D}v$ implies through Theorem \ref{BFT:Theorem:improved}-$(ii)$ that it suffices to find a constant $c(n,p)>0$ such that
\be\label{suffices_p<2}
 \mathcal{C}:=
  \inf_{v\in\test(B_1(0))\setminus\{0\}}
   \frac{\int_{B_1(0)}|x|^{p-n}|\nabla v|^pY^{2-p}_{k+1}(|x|/D) Y^{-1}_k(|x|/D)\dx}
          {\Big(\int_{B_1(0)}|x|^{-n}|v|^{p^\star}Y_k(|x|/D) X_{k+1}^{1+p^\star/p}(|x|/D)\dx\Big)^{p/p^\star}}
    \geq c(n,p).
\ee
\begin{remark}\label{rmrkGSHp<2}
Let $k=0$. Then the above sufficiency of \eqref{suffices_p<2} is straightforward from \eqref{GSHp<2} through the transform \eqref{GSH}; that is $u=f_{0,0}v$. It is for $k\in\N$ that we need the transform $u=f_{a,k,D}v$ for some $a>0$, hence Theorem \ref{BFT:Theorem:improved}-$(ii)$.
\end{remark}
To carry on with the proof, consider the Emden-Fowler transformation
\[
 v(x)=w(\tau,\theta),~~~~\mbox{where}~~\tau:=\frac{1}{X_{k+1}(r/D)},~~\theta:=\frac{x}{r}~~~~\mbox{with}~~r:=|x|.
\]
A simple computation using Lemma \ref{lemma:differential:rules} gives
\[
 \frac{\dd\tau}{\dd r}=-\frac{1}{r}Y_k(r/D),
\]
therefore,
\be\nonumber
 |\nabla v|^2
  =(\prt_rv)^2+\frac{1}{r^2}|\nabla_\theta v|^2
   =\frac{1}{r^2}Y_k^2(r/D)\big((\prt_\tau w)^2+Y_k^{-2}(r/D)|\nabla_{\theta}w|^2\big).
\ee
Let $F_1(t)$ denote the inverse function of $X_1(t)$ and define $F_{i+1}(t):=F_1\big(F_i(t)\big)$, $i=1, . . . , k$. With this notation, from the transformation we readily get
\[
 r/D= F_{k+1}(1/\tau),~~X_{i}(r/D)=F_{k+1-i}(1/\tau),~~i=1,...,k.
\]
Hence $Y_k(r/D)=\prod_{i=1}^{k}F_i(1/\tau)$. Setting $\tau_0:=X^{-1}_{k+1}(1/D)$, we deduce
\be\label{CKN_bc}
  \mathcal{C}
   =
    \inf_{\substack{w\in C^\infty([\tau_0,\infty)\times\us)\\w(\tau_0,\theta)=0}}
     \frac{\int_{\tau_0}^\infty\int_{\us}\tau^{p-2}\Big((\prt_\tau w)^2
             +\big(\prod_{i=1}^{k}F_i(1/\tau)\big)^{-2}|\nabla_{\theta}w|^2\Big)^{p/2}\dS(\theta)\dd\tau}
            {\Big(\int_{\tau_0}^{\infty}\int_{\us}\tau^{-1-p^*/p}|w|^{p^*}\dS(\theta)\dd \tau\Big)^{p/p^*}}.
 \ee

Suppose next that $n\geq3$ and set
\[
 \mathcal{S}
  :=
   \inf_{u\in\test(B_R)\setminus\{0\}}
    \frac{\int_{B_R}|\nabla u|^p~\dx}{\Big(\int_{B_R}|u|^{p^*}\dx\Big)^{p/p^*}}.
\]
From \cite{T} we know $\mathcal{S}=\mathcal{S}(n,p)>0$. Consider the transformation
\[
 u(x)=z(t,\theta),~~~~\mbox{where}~~t:=\frac{1}{r^{n-p}},~~\theta:=\frac{x}{r}~~~~\mbox{with}~~r:=|x|.
\]
An elementary computation gives
\[
 |\nabla u|^2=(\prt_r u)^2+\frac{1}{r^2}|\nabla_\theta u|^2 =(n-p)^2t^{2(n-p+1)/(n-p)}\Big((\prt_tz)^2+\big((n-p)t\big)^{-2}|\nabla_{\theta}z|^2\Big).
\]
Therefore,
\be\label{Sobolev_bc}
 \hspace{-1.8em}\frac{\mathcal{S}(n,p)}{(n-p)^{p(n-1)/n}}
  =
   \inf_{\substack{z\in C^\infty([R^{p-n},\infty)\times\us)\\z(R^{p-n},\theta)=0}}
    \frac{\int_{R^{p-n}}^{\infty}\int_{\us}t^{p-2}\Big((\prt_tz)^2+\big((n-p)t\big)^{-2}|\nabla_{\theta}z|^2\Big)^{p/2}\dS(\theta)\dd t}
         {\Big(\int_{R^{p-n}}^{\infty}\int_{\us}t^{-1-p^*/p}|z|^{p^*}\dS(\theta)\dd t\Big)^{p/p^*}}.
\ee
To compare the expressions on the right of \eqref{CKN_bc} and \eqref{Sobolev_bc}, we first choose $R$ such that $R^{p-n}=\tau_0$. Then we observe for $\tau\geq\tau_0$
\begin{align}\label{p<2, split}
 \tau^{-1}\prod_{i=1}^{k}F_i(1/\tau)
  & \leq \tau_0^{-1}\prod_{i=1}^{k}F_i(1/\tau_0)
    = \tau_0^{-1}Y_k(1/D)
     \leq 1
      \\ \nonumber &
       \leq n-p,
\end{align}
the last inequality because of $n\geq3$. Thus $\prod_{i=1}^{k}F_i(1/\tau)\leq(n-p)\tau$ for any $\tau\geq\tau_0$ and inserting this to \eqref{CKN_bc} we conclude with
\[
 \mathcal{C}\geq\frac{\mathcal{S}(n,p)}{(n-p)^{p(n-1)/n}}.
\]
This is \eqref{suffices_p<2} for $n\geq3$.

\smallskip

If $n=2$ we set
\[
 \mathcal{S}'
  :=\inf_{u\in C_0^{\infty}(B_R(0))\setminus\{0\}}
    \frac{\int_{B_R(0)}|x|^{\alpha p}|\nabla u|^p~\dx}
           {\Big(\int_{B_R(0)}|x|^{\alpha  p^\star}|u|^{p^\star}\dx\Big)^{p/p^\star}}.
\]
From \cite{CKN} (with $n=2$, $a=1$, $r=p^\star$ there), we know that
\[
 \mathcal{S}'=\mathcal{S}'(\alpha,p)>0~~~~\mbox{ whenever }~\al>1-2/p.
\]
In particular, taking $\alpha=1-1/p$ and considering the transformation
\[
 u(x)=z(t,\theta),~~~~\mbox{where}~~t:=\frac{1}{r},~~\theta:=\frac{x}{r}~~~~\mbox{with}~~r:=|x|,
\]
we deduce by a straightforward calculation
\be\label{Sobolev_bc_n=2}
 \mathcal{S}'(1-1/p,p)
   =
    \inf_{\substack{z\in C^\infty([R^{-1},\infty)\times\mathbb{S}^1)\\z(R^{-1},\theta)=0}}
     \frac{\int_{R^{-1}}^\infty \int_{\mathbb{S}^1}t^{p-2}\big((\prt_tz)^2+t^{-2}|\nabla_\theta z|^2\big)^{p/2}\dd\mathcal{H}^1(\theta)\dd t}
            {\Big(\int_{R^{-1}}^{\infty}\int_{\mathbb{S}^1}t^{-1-p^*/p}|z|^{p^*}\dd\mathcal{H}^1(\theta)\dd t\Big)^{p/p^*}}.
 \ee
To compare the expressions on the right of \eqref{CKN_bc} and \eqref{Sobolev_bc_n=2}, we choose $R$ such that $R^{-1}=\tau_0$. Then \eqref{p<2, split} says $\prod_{i=1}^{k}F_i(1/\tau)\leq\tau$ for any $\tau\geq\tau_0$  and inserting this to \eqref{CKN_bc} we conclude with
\[
 \mathcal{C}\geq\mathcal{S}'(1-1/p,p);
\]
that is \eqref{suffices_p<2} for $n=2$.\\

It remains to show the exponent $1+p^\star/p$ on $X_{k+1}$, $k\in\N$, cannot be decreased. The argument applies for any $1<p<n$. Suppose for the sake of contradiction, that $\e\in[0,1)$ is such that the following inequality holds for some $D\geq D_0$ and $k\in\N$
\be\label{contradiction1<p<2&kneq0}
 I_k[u]
  \geq
   c\Big(\int_{\Om}|u|^{p^\star}Y_k^{1+p^\star/p}(|x|/D) X_{k+1}^{(1+p^\star/p)\e}(|x|/D)\dx\Big)^{p/p^\star}~~~~\forall~u\in\test(\Om),
\ee
with $c$ being a positive constant independent of $u$. Applying H\"{o}lder's inequality with conjugate exponents $n/p$ and $p^\star/p$ and using (\ref{contradiction1<p<2&kneq0}) we have  (in the first displayed line below we write $Y_k$, $X_{k+1}$ instead of $Y_k(|x|/D)$, $X_{k+1}(|x|/D)$)
\begin{align}\nonumber
 \int_{\Om}\frac{|u|^p}{|x|^p}Y_k^2X_{k+1}^\gamma\dx
  & = \int_{\Om}\big\{|x|^{-p}Y_k^{1-p/p^\star}X_{k+1}^{\gamma-(1+p/p^\star)\e }\big\}\big\{|u|^pY_k^{1+p/p^\star}X_{k+1}^{(1+p/p^\star)\e }\big\}\dx
   \\ \label{contradictionIHI&kneq0} & \leq
    c^{-1}\Big(\int_{\Om}|x|^{-n}Y_k(|x|/D) X^\beta(|x|/D)\dx\Big)^{p/n}I_k[u],
\end{align}
where
\[
 \beta:=\Big[\gamma-\Big(1+\frac{p}{p^\star}\Big)\e\Big]\frac{n}{p}.
\]
The integral on the right is a constant depending on $n,p,\e,\gamma$ and $\Om$ if and only if $\beta>1$ (see for instance \cite[eq. (3.8)]{BFT2}). This is easily seen to be equivalent with
\be\label{betagamma}
 \gamma>2-\Big(1+\frac{p}{p^\star}\Big)(1-\e).
\ee
Thus for values of $\gamma$ determined from \eqref{betagamma}, we get from \eqref{contradictionIHI&kneq0} that
\[
 \int_{\Om}\frac{|u|^p}{|x|^p}Y_k^2(|x|/D)X_{k+1}^\gamma(|x|/D)\dx
  \leq
   c^{-1}C(n,p,\e,\gamma,\Om)I_k[u]~~~~\forall~u\in\test(\Om).
\]
However, Proposition 3.1-(i) of \cite{BFT1} (for the case $\kappa=n$ there) asserts the last inequality is possible only if $\gamma\geq2$. This and \eqref{betagamma} forces $\e\geq1$, a contradiction. \qed
\section{Theorem A for $p>2$}
We need the following special improvement to the series expansion of the $L^2$-Hardy inequality, which is valid only for radially symmetric functions.
\begin{lemma}\label{lemma:radial improvement}
Let $2\leq p<n$. Then for  any $D\geq1$, any $k\in\N\cup\{0\}$ and all radially symmetric functions $\zeta\in H_0^1\big(B_1(0)\big)$ we have
\begin{align*}\nonumber
 \int_{B_1(0)}|\nabla \zeta|^2\dx - \Big(\frac{n-2}{2}\Big)^2\int_{\Om}\frac{|\zeta|^2}{|x|^2}\dx
   - & \frac14\int_{B_1(0)}\frac{|\zeta|^2}{|x|^2}\sum_{i=1}^kY_i^2(|x|/D)\dx
   \\ & \geq C(n,p)\Big(\int_{B_1(0)}|x|^{p^\star(p-2)/p}Y_{k+1}^{1+p^\star/p}(|x|/D)|\zeta|^{2p^\star/p}~\dx\Big)^{p/p^\star}.
\end{align*}
\end{lemma}

\noindent\textbf{Proof.} We perform the change of variables
\be\label{GroundState_2_k}
 \zeta(r)=g_{k,D}(r)w(r),~~~~\mbox{where}~~g_{k,D}(r):=r^{1-n/2}Y^{-1/2}_k(r/D),~~r=|x|.
\ee
Then by Theorem \ref{BFT:Theorem:improved} for $p=2$ there, it is enough to prove that (it is in fact equivalent by Remark \ref{remark_equality})
\[
 \int_{B_1(0)}|x|^{2-n}Y_k^{-1}(|x|/D)|\nabla w|^2\dx
  \geq
   c(n,p)\Big(\int_{B_1(0)}|x|^{-n}Y_k(|x|/D)X_{k+1}^{1+p^\star/p}(|x|/D)|w|^{2p^\star/p}~\dx\Big)^{p/p^\star}.
\]
Since $w$ is radially symmetric, the above inequality is equivalent to
\[
 \int_0^1rY_k^{-1}(r/D)\big(w'(r)\big)^2\mathrm{d}r
  \geq
   c(n,p)\Big(\int_0^1r^{-1}Y_k(r/D)X_{k+1}^{1+p^\star/p}(r/D)|w(r)|^{2p^\star/p}~\mathrm{d}r\Big)^{p/p^\star}.
\]
The proof of this readily follows from \cite[Lemma 7.1]{FT} for $q=2p^\star/p$ there, or by \cite[Theorem 3 - pg 47]{Mz} for $\mathrm{d}\nu=rY_k^{-1}(r/D)\chi_{(0,1)}\mathrm{d}r$ and $\mathrm{d}\mu=r^{-1}Y_k(r/D)X_{k+1}^{1+p^\star/p}(r/D)\chi_{(0,1)}\mathrm{d}r$ there. \qed\\

\noindent{\bf{Proof of Theorem A for $p>2$.}} As in the case $p\leq2$, we can assume $\Om=B_1(0)$. Applying the transformation $u=f_{k,D}v$ we get by using both \eqref{BFTprime:p>2} and \eqref{BFTdoubleprime:p>2:Sobolev} that
\be\label{AB}
 I_k[u]
  \geq
   \frac{c_p}{2}\Big(\int_{\Om} f_{k,D}^p(x)|\nabla v|^p~\dx+\int_{\Om} g_{k,D}^2(x)|v|^{p-2}|\nabla v|^2\dx\Big),
\ee
where $g_{k,D}$ is given by \eqref{GroundState_2_k}. Note that $g_{k,D}$ is $f_{k,D}$ with $p=2$.

Following \cite{VzZ}, we use spherical coordinates $x=(r,\theta)$ ($r=|x|$ and $\theta=x/|x|$) to decompose $v(x)$ into spherical harmonics. For this purpose, let $\{h_l\}_{l\in\N\cup\{0\}}$ be the orthonormal basis of $L^2(\us)$ that is comprised of eigenfunctions of the Laplace-Beltrami operator $-\Delta_{\us}$ (the angular part of the Laplacian when expressed in spherical coordinates). This has corresponding eigenvalues $\lambda_l=l(l+n-2)$, $l\in\N\cup\{0\}$ (see \cite[Appendix]{Schn}). Thus
\[
 -\Delta_{\us}h_l=\lambda_lh_l~~\mbox{on }\us,
  ~~~~\mbox{and}~~~~
   \frac{1}{n\omn}\int_{\us}h_l(\theta)h_m(\theta)\dS(\theta)=\de_{lm}~~~~\mbox{for all }l,m\in\N\cup\{0\}.
\]
With these definitions we have the decomposition of $v\in\test(B_1(0))$ in its spherical harmonics
\[
 v(x)=\sum_{l=0}^\infty v_l(r)h_l(\theta).
\]
In particular $h_0(\theta)=1$ and the first term in the above decomposition is given by the spherical mean of $v$ on $\prt B_r(0)$, that is
\[
 v_0(r)=\frac{1}{n\omn r^{n-1}}\int_{\prt B_r(0)}v(x)\dS(x)=\frac{1}{n\omn}\int_{\us}v(r\theta)\dS(\theta).
\]
We now estimate the first term on the right hand side of \eqref{AB}. We have
\begin{align}\nonumber
 & \int_{\Om} f_{k,D}^p(x)|\nabla v|^p~\dx
  \\ \nonumber & = \int_0^1f_{k,D}^p(r)r^{n-1}\int_{\us}\Big((\prt_rv)^2+\frac{1}{r^2}|\nabla_\theta v|^2\Big)^{p/2}\dS(\theta)\mathrm{d}r
   \\ \label{A1A2} & \geq
    \int_0^1f_{k,D}^p(r)r^{n-1}\int_{\us}|\prt_rv|^p\dS(\theta)\mathrm{d}r
     +\underbrace{\int_0^1f_{k,D}^p(r)r^{n-1}\int_{\us}\frac{1}{r^p}|\nabla_\theta v|^p\dS(\theta)\mathrm{d}r}_{=:J},
\end{align}
by the fact that $(\kappa+\lambda)^q\geq \kappa^q+\lambda^q$, for all $\kappa,\lambda\geq0$ and any $q\geq1$. To estimate the first term on the right of \eqref{A1A2} we use \eqref{vecineq} to get
\begin{align}\nonumber
 \int_{\us}|\prt_rv|^p~\dS(\theta)
  & \geq \int_{\us}|\prt_rv_0|^p~\dS(\theta) + c_p\int_{\us}|\prt_r(v-v_0)|^p~\dS(\theta)
   \\\label{A11} & ~~~~ + p\int_{\us}|\prt_rv_0|^{p-2}(\prt_r v_0)\prt_r(v-v_0)\dS(\theta).
\end{align}
But since $\{v_l\}_{l\in\N\cup\{0\}}$ are radial
\begin{align*}
\int_{\us}|\prt_rv_0|^{p-2}(\prt_r v_0)\prt_r(v-v_0)\dS(\theta)
 & = |v'_0(r)|^{p-2}v'_0(r)\int_{\us}\partial_r(v-v_0)\dS(\theta)\\
  & = |v'_0(r)|^{p-2}v'_0(r)\sum_{l=1}^\infty v'_l(r)\int_{\us}f_l(\theta)\dS(\theta)
    = 0,
\end{align*}
and so
\[
 \int_{\us}|\prt_rv|^p~\dS(\theta)
  \geq
   c_p\int_{\us}|\prt_r(v-v_0)|^p~\dS(\theta),
\]
where we have cancel also the first term on the right hand side of \eqref{A11}. Plugging this to \eqref{A1A2} we deduce
\begin{align}\nonumber
 & \int_{\Om} f_{k,D}^p(x)|\nabla v|^p~\dx
  \\ \nonumber & \geq
   c_p\int_0^1f_{k,D}^p(r)r^{n-1}\int_{\us}\Big(|\prt_r(v-v_0)|^p+\frac{1}{r^p}|\nabla_\theta v|^p\Big) \dS(\theta)\mathrm{d}r+(1-c_p)J
  \\ \nonumber & \geq
   2^{1-p/2}c_p\int_0^1f_{k,D}^p(r)r^{n-1}\int_{\us}\Big(|\prt_r(v-v_0)|^2+\frac{1}{r^2}|\nabla_\theta v|^2\Big)^{p/2} \dS(\theta)\mathrm{d}r+(1-c_p)J
  \\ \label{A1A2A3} & =
   2^{1-p/2}c_p\int_{\Om} f_{k,D}^p(x)|\nabla(v-v_0)|^p~\dx+(1-c_p)J,
\end{align}
by the fact that $\kappa^q+\lambda^q\geq2^{1-q}(\kappa+\lambda)^q$, for all $\kappa,\lambda\geq0$ and any $q\geq1$. To estimate $J$ observe first that

\[
 \int_{\us}(v-v_0)\dS(\theta)=\sum_{l=1}^\infty v_l(r)\int_{\us}f_l(\theta)\dS(\theta)=0,
\]
so that utilizing once more the fact that $v_0$ is radial, we may use the Poincar\'e inequality on $\us$ (see for example \cite[Theorem 2.10]{H})
\begin{align}\nonumber
 \int_{\us}|\nabla_\theta v|^p~\dS(\theta)
  & = \int_{\us}|\nabla_\theta(v-v_0)|^p~\dS(\theta)
  \\ \nonumber & \geq
   C_{\mathbf{P}}(n,p)\int_{\us}|v-v_0|^p~\dS(\theta).
\end{align}
Inserting this in the definition of $J$, we get from \eqref{A1A2A3} the existence of a positive constant $C=C(n,p)$ such that
\begin{align}\nonumber
 \int_{\Om} f_{k,D}^p(x)|\nabla v|^p~\dx
  & \label{second term} \geq C(n,p)\Big(\int_{\Om} f_{k,D}^p(x)|\nabla (v-v_0)|^p~\dx
   +
    \int_{\Om} f_{k,D}^p(x)\frac{|v-v_0|^p}{|x|^p}\dx\Big)
  \\ & \geq C(n,p)\Big(\int_{\Om}|x|^{-n}Y_k(|x|/D) X_{k+1}^{1+p^\star/p}(|x|/D)|v-v_0|^{p^\star}~\dx\Big)^{p/p^\star},
\end{align}
where in the last inequality we have used the Sobolev inequality and the fact that $X_i\leq1$ for all $i\in\N$.

\medskip

Next we estimate the second term on the right hand side of \eqref{AB}. Setting $w=|v|^{p/2}$ we have
\begin{align}
\int_{\Om} g_{k,D}^2(x)|v|^{p-2}|\nabla v|^2\dx=\int_{\Om} g_{k,D}^2(x)|\nabla w|^2\dx.\label{v to w}
\end{align}
Now we assert that the function $\zeta=g_{k,D}w$ belongs to $H_0^1(\Om)$. Indeed
\[
 \zeta=g_{k,D}w=g_{k,D}|v|^{p/2}=g_{k,D}f_{k,D}^{-p/2}|u|^{p/2}=|x|^{1-p/2}|u|^{p/2}\in H_0^1(\Om),
 \]
since $2\leq p<n$. Thus by Remark \ref{remark_equality} we have that the following equality is valid
\be\label{w to l2series}
 \int_{\Om} g_{k,D}^2(x)|\nabla w|^2\dx
  =
  \int_{\Om}|\nabla \zeta|^2dx-\Big(\frac{n-2}{2}\Big)^2\int_{\Om}\frac{|\zeta|^2}{|x|^2}\dx
   -
    \frac14\int_{\Om}\frac{|\zeta|^2}{|x|^2}\sum_{i=1}^kY_i^2(|x|/D)\dx.
\ee
Taking the decomposition of $\zeta$ in its spherical harmonics we know that (see \cite[eq. (7.6)]{FT})
\begin{align}\nonumber
 & \int_{\Om}|\nabla \zeta|^2\dx-\Big(\frac{n-2}{2}\Big)^2\int_{\Om}\frac{|\zeta|^2}{|x|^2}\dx
  - \frac14\int_{\Om}\frac{|\zeta|^2}{|x|^2}\sum_{i=1}^kY_i^2(|x|/D)\dx
   \\ \nonumber & \geq \int_{\Om}|\nabla \zeta_0|^2\dx - \Big(\frac{n-2}{2}\Big)^2\int_{\Om}\frac{|\zeta_0|^2}{|x|^2}\dx
      - \frac14\int_{\Om}\frac{|\zeta_0|^2}{|x|^2}\sum_{i=1}^kY_i^2(|x|/D)\dx
      \\ \label{L2series to radSob} & \geq C(n,p) \Big(\int_{\Om}|x|^{p^\star(p-2)/p}Y_{k+1}^{1+p^\star/p}(|x|/D)|\zeta_0|^{2p^\star/p}\;\dx\Big)^{p/p^\star},
\end{align}
where in the last inequality we have used Lemma \ref{lemma:radial improvement} since $\zeta_0$ is radial. In particular, we have that
\[
 \zeta_0(r)=\frac{1}{n\omn r^{n-1}}\int_{\partial B_r(0)}\zeta(x)\dS(x),
\]
which implies (note that $\zeta_0$ is nonnegative since $\zeta$ is nonnegative)
\begin{align}\nonumber
 \zeta_0(r)
  & = \frac{g_{k,D}(r)}{n\omn r^{n-1}}\int_{\prt B_r(0)}w(x)\dS(x)
   \\ \nonumber & = \frac{g_{k,D}(r)}{n\omn r^{n-1}}\int_{\partial B_r(0)}|v(x)|^{p/2}\dS(x)
    \\ \nonumber & \geq g_{k,D}(r)\Big|\frac{1}{n\omn r^{n-1}}\int_{\prt B_r(0)}v(x)\dS(x)\Big|^{p/2}
     \\ \nonumber & = g_{k,D} |v_0(r)|^{p/2}.
\end{align}
This applied to \eqref{L2series to radSob} gives together with \eqref{w to l2series} and \eqref{v to w}  that
\be\label{first term}
 \int_{\Om} g_{k,D}^2(x)|v|^{p-2}|\nabla v|^2\dx
  \geq
   c(n,p)\Big(\int_{\Om}|x|^{-n}Y_k(|x|/D) X_{k+1}^{1+p^\star/p}(|x|/D)|v_0|^{p^\star}\dx\Big)^{p/p^*}.
\ee
Inserting \eqref{first term} and \eqref{second term} in \eqref{AB}, inequality \eqref{Hardy-Sobolev-k:L^p<n} follows. Moreover, it is proved in the previous section that the exponent $1+p^\star/p$ on $X_{k+1}$, $k\in\N$, cannot be decreased. \qed
\section{A local estimate}
The local estimate of Theorem \ref{thrm:local} below is the key estimate in order to establish the series improvement to the Hardy-Morrey inequality that appears in Theorem B. To establish it we need the following weighted Hardy inequality with trace term.

\begin{lemma}\label{Lemma:trace}
Let $\gamma\in\rr\setminus\{0\}$ and $U$ be a bounded domain in $\rn$, $n\geq2$, having locally Lipschitz boundary. Denote by $\nu(x)$ the exterior unit normal vector defined at almost every $x\in\prt U$. Then for all $D\geq R_U:=\sup_{x\in U}|x|$, $q\geq1$, $k\in\N$, all $s\neq n$ and any $v\in\test(\rn\setminus\{0\})$, there holds
\begin{align}\nonumber
 \Big|\frac{q}{n-s}\Big|^q\int_U\frac{|\nabla v(x)|^q}{|x|^{s-q}}Y_k^\gamma(|x|/D)\dx
  & +
   \frac{q}{n-s}
    \int_{\prt U}\frac{|v(x)|^q}{|x|^s}Y_k^\gamma(|x|/D) x\cdot\nu(x)\dS(x)
     \\ \label{weighted:Hardy:trace} & ~~~~ \geq
       \int_U\frac{|v(x)|^q}{|x|^s}Y_k^\gamma(|x|/D) \Big[1+\frac{\gamma q}{n-s}Z_k(|x|/D)\Big]\dx.
\end{align}
\end{lemma}
\noindent\textbf{Proof.} Integration by parts together with Lemma \ref{lemma:differential:rules} give
\begin{align}\nonumber &
 -\int_U|v|Y_k^\gamma(|x|/D)\Div\Big\{\frac{x}{|x|^s}\Big\}\dx
   \\ \nonumber & =
    \int_U\nabla|v|\cdot\frac{x}{|x|^s}Y_k^\gamma(|x|/D)\dx
     +\gamma\int_U\frac{|v|}{|x|^s}Y_k^\gamma(|x|/D) Z_k(|x|/D)\dx
      -\int_{\prt U}|v|Y_k^\gamma(|x|/D) \frac{x}{|x|^s}\cdot\nu\dS(x),
\end{align}
and since $\Div\{|x|^{-s}x\}=(n-s)|x|^{-s}$, $x\neq0$, we get
\[
 \int_U\frac{|\nabla v|}{|x|^{s-1}}Y_k^\gamma(|x|/D)\dx
  -
   \int_{\prt U}\frac{|v|}{|x|^s}Y_k^\gamma(|x|/D)  x\cdot\nu\dS(x)
    \geq
     \int_U\frac{|v|}{|x|^s}Y_k^\gamma(|x|/D)\big[s-n-\gamma Z_k(|x|/D)\big]\dx,
\]
if $s>n$, or
\[
 \int_U\frac{|\nabla v|}{|x|^{s-1}}Y_k^\gamma(|x|/D)\dx
  +
   \int_{\prt U}\frac{|v|}{|x|^s}Y_k^\gamma(|x|/D) x\cdot\nu\dS(x)
    \geq
     \int_U\frac{|v|}{|x|^s}Y_k^\gamma(|x|/D)\big[n-s+\gamma Z_k(|x|/D)\big]\dx,
\]
if $s<n$, where we have also used the fact that $|\nabla|v(x)||\leq|\nabla v(x)|$ for a.e. $x\in U$. We may write both inequalities in one as follows
\begin{align}\nonumber
 \frac{1}{|n-s|}\int_U\frac{|\nabla v|}{|x|^{s-1}}Y_k^\gamma(|x|/D)\dx
  & +
    \frac{1}{n-s}\int_{\prt U}\frac{|v|}{|x|^s}Y_k^\gamma(|x|/D) x\cdot\nu\dS(x)
    \\ \nonumber & ~~~~ \geq
     \int_U\frac{|v|}{|x|^s}Y_k^\gamma(|x|/D) \Big[1+\frac{\gamma}{n-s}Z_k(|x|/D)\Big]\dx.
\end{align}
This is inequality \eqref{weighted:Hardy:trace} for $q=1$. Substituting $v$ by $|v|^q$ with $q>1$, we arrive at
 \begin{align}\nonumber
  \frac{q}{|n-s|}\int_U\frac{|\nabla v||v|^{q-1}}{|x|^{s-1}}Y_k^\gamma(|x|/D)\dx
   & +
    \frac{1}{n-s}
     \int_{\prt U}\frac{|v|^q}{|x|^s}Y_k^\gamma(|x|/D)  x\cdot\nu\dS(x)
      \\ \label{Lemma:trace:1-q} & ~~~~ \geq
       \int_U\frac{|v|^q}{|x|^s}Y_k^\gamma(|x|/D) \Big[1+\frac{\gamma}{n-s}Z_k(|x|/D) \Big]\dx.
 \end{align}
The first term on the left of \eqref{Lemma:trace:1-q} can be written as follows
 \begin{align}\nonumber
  \frac{q}{|n-s|}\int_U\frac{|\nabla v||v|^{q-1}}{|x|^{s-1}}Y_k^\gamma(|x|/D)\dx
   & =
    \int_U\Big\{\frac{q}{|n-s|}\frac{|\nabla v|}{|x|^{s/q-1}}\Big\}
     \Big\{\frac{|v|^{q-1}}{|x|^{s-s/q}}\Big\}Y_k^\gamma(|x|/D)\dx
      \\ \nonumber & \hspace{-5em} \leq
       \frac{1}{q}\Big|\frac{q}{n-s}\Big|^{q}
        \int_U\frac{|\nabla v|^q}{|x|^{s-q}}Y_k^\gamma(|x|/D)\dx
         +
          \frac{q-1}{q}\int_U\frac{|v|^q}{|x|^s}Y_k^\gamma(|x|/D)\dx,
 \end{align}
by Young's inequality. Thus \eqref{Lemma:trace:1-q} becomes
\begin{align}\nonumber
 \frac{1}{q}\Big|\frac{q}{n-s}\Big|^{q}
  \int_U\frac{|\nabla v|^q}{|x|^{s-q}}Y_k^\gamma(|x|/D)\dx
   & +
    \frac{1}{n-s}
     \int_{\prt U}\frac{|v|^q}{|x|^s}Y_k^\gamma(|x|/D)  x\cdot\nu\dS(x)
      \\ \nonumber & ~~~~ \geq
       \int_U\frac{|v|^q}{|x|^s}Y_k^\gamma(|x|/D) \Big[\frac{1}{q}+\frac{\gamma}{n-s}Z_k(|x|/D) \Big]\dx.
\end{align}
Multiplying by $q$ we get \eqref{weighted:Hardy:trace}.\qed

\begin{theorem}\label{thrm:local} Let $2\leq p\neq n$ and $1\leq q<p$. There exist constants $B=B(n,p,q)\geq1$ and $C=C(n,p,q)>0$ such that for all $u\in\test(\Om\setminus\{0\})$, any ball $B_r$ of radius $r\in\big(0,\diam(\Om)\big)$ that contains the origin, any $D\geq B\diam(\Om)$ and all $k\in\N$
 \begin{align}\label{local}
  \int_{B_r}\frac{|u|^q}{|x|^q}\Big[1-\frac{q^2}{n(p-q)}Z_k(|x|/D)\Big]\dx
   \leq
    Cr^{n(1-q/p)}Y_{k+1}^{-q/p}(r/D)\big(I_k[u;D]\big)^{q/p}.
 \end{align}
\end{theorem}

\noindent\textbf{Proof.}  Given $u\in\test(\Om\setminus\{0\})$ we define as usual $v\in\test(\Om\setminus\{0\})$ through the transform \[u(x)=f_{k,D}(x)v(x),\]
where $D\geq\diam(\Om)$ and $1<p\neq n$. Then with $q\in[1,p)$ and $r\in(0,\diam(\Om))$, we have for any ball containing the origin that
 \begin{align}\nonumber
  \int_{B_r}\frac{|u|^q}{|x|^q}\Big[1-\frac{q^2}{n(p-q)}Z_k(|x|/D)\Big]\dx
   & =
    \int_{B_r}\frac{|v|^q}{|x|^{nq/p}}Y_k^{-q/p}(|x|/D) \Big[1-\frac{q^2}{n(p-q)}Z_k(|x|/D)\Big]\dx
     \\ \nonumber & \leq
      \Big(\frac{pq}{n(p-q)}\Big)^q\underbrace{\int_{B_r}|x|^{q(p-n)/p}|\nabla v|^qY_k^{-q/p}(|x|/D)\dx}_{=:M_r}
       \\ \label{Estim:Mr:Pr} & ~~~~ + \frac{pq}{n(p-q)}
          \underbrace{\int_{\prt B_r}\frac{|v|^q}{|x|^{nq/p}}Y_k^{-q/p}(|x|/D)  x\cdot\nu\dS(x)}_{=:P_r},
 \end{align}
where we have used Lemma \ref{Lemma:trace} for $U=B_r$, $s=nq/p$ and $\gamma=-q/p$. By H\"{o}lder's inequality
 \begin{align}\nonumber
  M_r
   & \leq
    \big(\omn r^n\big)^{1-q/p}
     \Big(\int_{B_r}|x|^{p-n}|\nabla v|^pY_k^{-1}(|x|/D)\dx\Big)^{q/p}
      \\ \nonumber \mbox{(by (\ref{BFTprime:p>2}))}~~ & \leq
       C(n,p,q)r^{n(1-q/p)}\big(I_k[u;D]\big)^{q/p}
        \\ \label{Estim:Mr} & \leq
         C(n,p,q)r^{n(1-q/p)}Y_{k+1}^{-q/p}(r/D)\big(I_k[u;D]\big)^{q/p},
 \end{align}
for any $D\geq b'\diam(\Om)$, $b'=b'(n,p)\geq1$. The last inequality is true since $0<Y_{k+1}(t)\leq1$ for all $t\in(0,1]$. For $\it{P}_r$, noting that $x\cdot\nu\geq0$ for all $x\in\prt B_r$ ($B_r$ is star-shaped with respect to any of it's points; thus $0$ particular), we may also apply H\"{o}lder's inequality as follows
 \begin{align}\nonumber &
  \textit{P}_r
   =
    \int_{\prt B_r}\big\{Y_{k+1}^{-q/p}(|x|/D)\big\}
     \Big\{\frac{|v|^q}{|x|^{nq/p}}X_{k+1}^{q/p}(|x|/D) \Big\}x\cdot\nu\dS(x)
       \\ \label{Estim:Sr:Tr} & \hspace{-3em}\leq
        \Big(\underbrace{\int_{\prt B_r}Y_{k+1}^{-q/(p-q)}(|x|/D)  x\cdot\nu\dS(x)}_{=:S_r}\Big)^{1-q/p}
         \Big(\underbrace{\int_{\prt B_r}\frac{|v|^p}{|x|^n}X_{k+1}(|x|/D)  x\cdot\nu\dS(x)}_{=:T_r}\Big)^{q/p}.
 \end{align}
By the divergence theorem we have
 \begin{align}\nonumber
  S_r
   & =
    \int_{B_r}\Div\Big\{Y_{k+1}^{-q/(p-q)}(|x|/D) ~x\Big\}\dx
     \\ \nonumber & =
      n\int_{B_r}Y_{k+1}^{-q/(p-q)}(|x|/D)\dx
       -\frac{q}{p-q}\int_{B_r}Y_{k+1}^{-q/(p-q)}(|x|/D)  Z_k(|x|/D)\dx
        \\ \nonumber & \leq
         n\int_{B_r(0)}Y_{k+1}^{-q/(p-q)}(|x|/D)\dx,
 \end{align}
since this integral increases if we change the domain of integration from $B_r$ to $B_r(0)$. Thus
 \begin{align}\nonumber
  S_r
   & \leq
    n^2\omn\int_0^rt^{n-1}Y_{k+1}^{-q/(p-q)}(t/D)\mathrm{d}t
     \\ \label{estim:Sr} & \leq
      C(n)r^nY_{k+1}^{-q/(p-q)}(r/D),
 \end{align}
for any $D\geq {\eta}\diam(\Om)$, $\eta\geq1$ depending only on $n,p,q$, by a direct application of Lemma \ref{Lemma:X} for $\al=n$ and $\beta=q/(p-q)$. To estimate $\it{T}_r$ we also employ the divergence theorem to get
 \[
  T_r
   =
      \int_{B_r}\Div\big\{|x|^{-n}X_{k+1}(|x|/D)  x\big\}|v|^p\dx
      +\int_{B_r}|x|^{-n}X_{k+1}(|x|/D)  x\cdot\nabla(|v|^p)\dx.
 \]
A direct computation using Lemma \ref{lemma:differential:rules} shows that
 \[
  \Div\big\{|x|^{-n}X_{k+1}(|x|/D)x\big\}
   =
    |x|^{-n}Y_k(|x|/D)X^2_{k+1}(|x|/D),
     ~~~~x\in\Om\setminus\{0\}.
 \]
Returning then to the original function $u$ in the first integral and taking the absolute value in the second, we arrive at
 \begin{align}\label{Estim:Tr:1}
  T_r
    \leq
    \int_{\Om}\frac{|u|^p}{|x|^p}Y^2_{k+1}(|x|/D)\dx
     +p\underbrace{\int_{\Om}\frac{|v|^{p-1}}{|x|^{n-1}}|\nabla v|X_{k+1}(|x|/D)\dx}_{=:J}.
 \end{align}
Now we apply the Cauchy-Schwarz inequality in the second integral above as follows
 \begin{align}\nonumber
  J
   & =
    \int_{\Om}\Big\{\frac{|v|^{p/2-1}}{|x|^{n/2-1}}|\nabla v|Y_k^{-1/2}(|x|/D) \Big\}
      \Big\{\frac{|v|^{p/2}}{|x|^{n/2}}Y_k^{1/2}(|x|/D) X_{k+1}(|x|/D) \Big\}\dx
       \\ \nonumber & \leq
        \Big(\int_{\Om}\frac{|v|^{p-2}}{|x|^{n-2}}|\nabla v|^2Y_k^{-1}(|x|/D)\dx\Big)^{1/2}
          \Big(\int_{\Om}\frac{|v|^p}{|x|^n}Y_k(|x|/D) X^2_{k+1}(|x|/D)\dx\Big)^{1/2}
           \\ \nonumber & \leq
            c(p)\big(I_k[u;D]\big)^{1/2}
             \Big(\int_{\Om}\frac{|u|^p}{|x|^p}Y_{k+1}^2(|x|/D)\dx\Big)^{1/2},
\end{align}
for all $D\geq b'''D_0$, where $b'''\geq1$ depends only on $n,p$. Here we have used \eqref{BFTdoubleprime:p>2:Sobolev} to estimate the first factor and returned to the original function $u$ in the second factor. Estimate \eqref{Estim:Tr:1} is now
\be\nonumber
 T_r
  \leq
   \int_{\Om}\frac{|u|^p}{|x|^p}Y_{k+1}^2(|x|/D)\dx
    +c(p)\big(I_k[u;D]\big)^{1/2}
     \Big(\int_{\Om}\frac{|u|^p}{|x|^p}Y_{k+1}^2(|x|/D)\dx\Big)^{1/2}.
\ee
According to Theorem \ref{BFT:Theorem}, there exist constants $b=b(n,p)\geq1$ and $c(n,p)>0$, both depending only on $n,p$, such that for any $D\geq bD_0$, the common integral appearing on the right hand side is bounded above by $c(n,p)I_k[u;D]$. It follows that
 \be\label{Estim:Tr}
  T_r
   \leq
    C(n,p)I_k[u;D],
 \ee
for any $D\geq\max\{b,b'''\}D_0$. Setting $b''=\max\{b,b''',\eta\}$ and since $0\in\Om$ implies $D_0\leq\diam(\Om)$, we may insert \eqref{Estim:Tr} into \eqref{Estim:Sr:Tr} and taking into account \eqref{estim:Sr} we end up with
 \be\nonumber
  P_r
   \leq
    C(n,p,q)r^{n(1-q/p)}Y_{k+1}^{-q/p}(r/D)\big(I_k[u;D]\big)^{1/p},
 \ee
for any $D\geq b''\diam(\Om)$. The last inequality together with \eqref{Estim:Mr}, when applied to estimate \eqref{Estim:Mr:Pr} gives \eqref{local} for any $D\geq B\diam(\Om)$ with $B=\max\{b',b''\}$. \qed
\section{Proof of Theorem B}
We start with \eqref{k:improvedHardy:Morrey} when one point in the H\"{o}lder semi-norm taken to be the origin.

\begin{proposition} Let $p>n$. There exist constants $\tilde{B}=\tilde{B}(n,p)\geq1$ and $C=C(n,p)>0$ such that for any $k\in\N\cup\{0\}$, all $D\geq\tilde{B}\diam(\Om)$ and all $u\in\test(\Om\setminus\{0\})$
 \be\label{k:improvedHardy:Morrey:onepoint}
  \sup_{x\in\Om}\Big\{\frac{|u(x)|}{|x|^{1-n/p}}Y_{k+1}^{1/p}(|x|/D)\Big\}
   \leq C
    \big(I_k[u;D]\big)^{1/p}.
 \ee
\end{proposition}

\noindent\textbf{Proof.} Let $B_r$ be a ball of radius $r\in(0,\diam(\Om))$ and set
\[
 u_{B_r}:=\frac{1}{|B_r|}\int_{B_r}u(z)\dd z.
\]
By the local version of Sobolev's integral representation formula (see \cite{GTr}-Lemma 7.16), we have
 \[
  |u(x)-u_{B_r}|
    \leq \frac{2^n}{n\omn}
    \int_{B_r}\frac{|\nabla u(z)|}{|x-z|^{n-1}}\dd z,~~~~x\in B_r.
 \]
Applying the transform $u(z)=f_{k,D}(z)v(z)$, we get
\begin{align}\nonumber
  \frac{n\omn}{2^n} |u(x)-u_{B_r}|
   & \leq
    \int_{B_r}\frac{|z|^{1-n/p}Y_k^{-1/p}(|z|/D)|\nabla v(z)|}{|x-z|^{n-1}}\dd z
     +\int_{B_r}\frac{Y_k^{-1/p}(|z|/D)|A_{0,k}(|z|/D)||v(z)|}{|z|^{n/p}|x-z|^{n-1}}\dd z
     \\ \label{Estim:Kr:Lr} & =:
      K_r(x)+L_r(x),
 \end{align}
with $A_{0,k}$ given by \eqref{def:Ak} with $a=0$. By H\"{o}lder's inequality
\begin{align}\nonumber
  K_r(x)
   & \leq
    \Big(\int_{B_r}\frac{1}{|x-z|^{(n-1)p/(p-1)}}\dd z \Big)^{1-1/p}
     \Big(\int_{B_r}|z|^{p-n}Y_k^{-1}(|z|/D)|\nabla v|^p\dd z \Big)^{1/p}
 \\ \nonumber & \leq
  \Big(\int_{B_r(x)}\frac{1}{|x-z|^{(n-1)p/(p-1)}}\dd z \Big)^{1-1/p}
    \Big(\int_{\Om}|z|^{p-n}Y_k^{-1}(|z|/D)|\nabla v|^p\dd z \Big)^{1/p}.
\end{align}
Using now \eqref{BFTprime:p>2} we obtain the following estimate on $K_r$
 \begin{align}\nonumber
  K_r(x)
   & \leq C(n,p)
    r^{1-n/p}\big(I_k[u;D]\big)^{1/p}
     \\ \label{Estim:Kr} & \leq C(n,p)
      r^{1-n/p}Y_{k+1}^{-1/p}(r/D)\big(I_k[u;D]\big)^{1/p},~~~~x\in B_r,
 \end{align}
for any $D\geq\diam(\Om)$, where the last inequality is a consequence of $0<Y_{k+1}(t)\leq1$ for all $t\in(0,1]$. Next we fix $0<\e<(p-n)/n$ and estimate $L_r(x)$. By H\"{o}lder's inequality
 \be\label{Estim:Lr:1}
  L_r(x)
   \leq
    \Big(\underbrace{\int_{B_r}\frac{|A_{0,k}(|z|/D)|}{|x-z|^{(n-1)p/(p-1-\e)}}\dd z}_{=:M_{r,D}(x)}\Big)^{1-(1+\e)/p}
     \Big(\int_{B_r}\frac{|v|^{p/(1+\e)}}{|z|^{n/(1+\e)}}|A_{0,k}(|z|/D)|\dd z\Big)^{(1+\e)/p},~~~~x\in B_r.
 \ee
Assumption $\e<(p-n)/n$ guarantees $M_{r,D}(x)< \infty$ for all $x\in B_r$. More precisely, recalling first Remark \ref{rmrk_convergence}, we may restrict $D$ so that $D\geq b'''D_0$ with some $b'''=b'''(n,p)\geq1$ so that $A_{0,k}(|z|/D)\geq0$ for all $x\in B_r$. Then we have
\begin{align}\nonumber
  M_{r,D}(x)
   & \leq
    \frac{p-n}{p}\int_{B_r(x)}\frac{1}{|x-z|^{(n-1)p/(p-1-\e)}}\dd z
     \\ \nonumber & =
      C(n,p)r^{(p-n-n\e)/(p-1-\e)},~~~~x\in B_r.
 \end{align}
Returning to the original function $u$ on the right of \eqref{Estim:Lr:1}, we obtain for all $D\geq b'''D_0$ that
 \be\label{Estim:Lr:3/2}
  L_r(x)
   \leq C(n,p)
    r^{1-n/p-n\e/p}\Big(\int_{B_r}\frac{|u|^{p/(1+\e)}}{|z|^{p/(1+\e)}}A_{0,k}(|z|/D)\dd z\Big)^{(1+\e)/p},~~~~x\in B_r.
 \ee
At this point we use Theorem \ref{thrm:local} with $q=p/(1+\e)$; that is
\be\label{local-q/enakiepsilon}
 \int_{B_r}\frac{|u|^{p/(1+\e)}}{|x|^{p/(1+\e)}}\Big[1-\frac{p}{n\e(1+\e)}Z_k(|z|/D)\Big]\dx
  \leq
   C(n,p)r^{n\e/(1+\e)}Y_{k+1}^{-1/(1+\e)}(r/D)\big(I_k[u;D]\big)^{1/(1+\e)}.
\ee
To couple this with \eqref{Estim:Lr:3/2} we need a positive constant $\lambda=\lambda(n,p)$ such that
\[
 A_{0,k}(|z|/D)
  \leq
   \lambda\Big[1-\frac{p}{n\e(1+\e)}Z_k(|z|/D)\Big],~~~~\forall~z\in\Om.
\]
Taking any $\lambda$ such that $\lambda>(p-n)/p$, keeping in mind that $\e<(p-n)/n$ and recalling the definition of $A_{0,k}$, this is the same as
\be\label{zk-restriction}
 Z_k(|z|/D)
  \leq
   \frac{\lambda-\frac{p-n}{p}}{\lambda\frac{p}{n\e(1+\e)}-\frac{1}{p}},~~~~\forall~z\in\Om.
\ee
which is satisfied after a possible further restriction on $D$. More precisely, note once more that because of Remark \ref{remark:Zinfinite} we can achieve \eqref{zk-restriction} for sufficiently large $\bar{b}=\bar{b}(n,p)$ and all $D\geq\bar{b} D_0$. Plugging \eqref{local-q/enakiepsilon} to \eqref{Estim:Lr:3/2} we obtain
\[
  L_r(x)
   \leq C(n,p)
    r^{1-n/p-n\e/p}\Big(\int_{B_r}\frac{|u|^{p/(1+\e)}}{|z|^{p/(1+\e)}}\Big[1-\frac{p}{n\e(1+\e)}Z_k(|z|/D)\Big]\dd z\Big)^{(1+\e)/p},~~~~x\in B_r,
 \]
for all $D\geq\max\{b''',\bar{b}\}D_0$. Using Theorem \ref{thrm:local} with $q=p/(1+\e)$,
 \begin{align}\nonumber
  L_r(x)
   & \leq C(n,p)
    r^{1-n/p-n\e/p}\Big(r^{n\e/(1+\e)}Y_{k+1}^{-1/(1+\e)}(r/D)\big(I_k[u;D]\big)^{1/(1+\e)}\Big)^{(1+\e)/p}
     \\ \label{Estim:Lr:3} & = C(n,p)
      r^{1-n/p}Y_{k+1}^{-1/p}(r/D)\big(I_k[u;D]\big)^{1//p},
  \end{align}
for any $D\geq\tilde{B}\diam(\Om)$, where $\tilde{B}$ depends only on $n,p$.

Applying estimates (\ref{Estim:Lr:3}) and (\ref{Estim:Kr}) to estimate (\ref{Estim:Kr:Lr}), we conclude
 \be\label{Final}
  |u(x)-u_{B_r}|
   \leq C(n,p)
    r^{1-n/p}Y_{k+1}^{-1/p}(r/D)\big(I_k[u;D]\big)^{1/p},
 \ee
for all $x\in B_r$ and any $D\geq\tilde{B}\diam(\Om)$. Since $0\in B_r$, it follows from (\ref{Final}) that
 \[
  |u_{B_r}|
   \leq C(n,p)
    r^{1-n/p}Y_{k+1}^{-1/p}(r/D)\big(I_k[u;D]\big)^{1/p}.
 \]
Hence
 \begin{align}\nonumber
  |u(x)|
   & \leq
    |u(x)-u_{B_r}| + |u_{B_r}|
     \\ \nonumber & \leq
      C(n,p) r^{1-n/p}Y_{k+1}^{-1/p}(r/D)\big(I_k[u;D]\big)^{1/p},
 \end{align}
for all $x\in B_r$ and any $D\geq\tilde{B}\diam(\Om)$. Now if $x\in\Om$ we may consider a ball $B_r$ of radius $r=3|x|/2$ containing $x$ and the previous inequality yields \eqref{k:improvedHardy:Morrey:onepoint}. \qed \\

\noindent\textbf{Proof of Theorem B.} Let $x,y\in\Om$, $x\neq y$, and consider a ball $B_r$ of radius $r:=|x-y|$ that contains $x,y$. Then $r\in(0,\diam(\Om))$ and we have
 \begin{align}\nonumber
  |u(x)-u(y)|
   & \leq
    |u(x)-u_{B_r}| + |u(y)-u_{B_r}|
     \\ \label{Estim:J} & \leq \frac{2^n}{n\omn}
      \Big\{\underbrace{\int_{B_r}\frac{|\nabla u(z)|}{|x-z|^{n-1}}\dd z}_{=:J(x)}
       +\underbrace{\int_{B_r}\frac{|\nabla u(z)|}{|y-z|^{n-1}}\dd z}_{=:J(y)}\Big\},
 \end{align}
where we have used Sobolev's integral representation formula (see \cite{GTr}-Lemma 7.16) twice. In what follows we estimate $J(x)$ independently of $x$. Applying the transform $u(z)=f_{k,D}(z)v(z)$, we get
 \be\label{Estim:Krr:Lrr}
  J(x)
   \leq
    \underbrace{\int_{B_r}\frac{|z|^{1-n/p}Y_k^{-1/p}(|z|/D)|\nabla v(z)|}{|x-z|^{n-1}}\dd z}_{=:\mathrm{K}_r(x)}
     +\underbrace{\int_{B_r}\frac{Y_k^{-1/p}(|z|/D)|A_{0,k}(|z|/D)||v(z)|}{|z|^{n/p}|x-z|^{n-1}}\dd z}_{=:\Lambda_r(x)}.
 \ee
Working as we did to get \eqref{Estim:Kr} we obtain
\be\label{Estim:Krr}
  \mathrm{K}_r(x)
   \leq C_1(n,p)
    r^{1-n/p}Y_{k+1}^{-1/p}(r/D)\big(I_k[u;D]\big)^{1/p},
 \ee
for any $D\geq\diam(\Om)$. Next we rewrite $\Lambda_r(x)$ with the original function $u$ to get
 \[
  \Lambda_r(x)
   =
    \int_{B_r}\frac{|A_{0,k}(|z|/D)||u(z)|}{|z||x-z|^{n-1}}\dd z.
 \]
We insert \eqref{k:improvedHardy:Morrey:onepoint} in $\Lambda_r(x)$ to deduce
 \[
  \Lambda_r(x)
   \leq C(n,p)\big(I_k[u;D]\big)^{1/p}
    \int_{B_r}\frac{Y_{k+1}^{-1/p}(|z|/D)|A_{0,k}(|z|/D)|}{|z|^{n/p}|x-z|^{n-1}}\dd z,
 \]
for any $D\geq\tilde{B}\diam(\Om)$. Recalling once more Remark \ref{rmrk_convergence}, we can further restrict $D$ so that $D\geq \max\{b''',\tilde{B}\}\diam(\Om)$ and then
 \[
  \Lambda_r(x)
   \leq C(n,p)\big(I_k[u;D]\big)^{1/p}
    \int_{B_r}\frac{Y_{k+1}^{-1/p}(|z|/D)}{|z|^{n/p}|x-z|^{n-1}}\dd z,
 \]
Letting $n<Q<p$ we may use H\"{o}lder's inequality to obtain
 \begin{align}\nonumber
  \Lambda_r(x)
   & \leq C(n,p)\big(I_k[u;D]\big)^{1/p}
    \Big(\int_{B_r}\frac{Y_{k+1}^{-Q/p}(|z|/D)}{|z|^{nQ/p}}\dd z\Big)^{1/Q}
     \Big(\int_{B_r}\frac{1}{|x-z|^{(n-1)Q'}}\dd z\Big)^{1/Q'}
      \\ \nonumber & \leq
        C(n,p)\big(I_k[u;D]\big)^{1/p}
         \Big(\int_{B_r(0)}\frac{Y_{k+1}^{-Q/p}(|z|/D)}{|z|^{nQ/p}}\dd z\Big)^{1/Q}
           \Big(\int_{B_r(x)}\frac{1}{|x-z|^{(n-1)Q'}}\dd z\Big)^{1/Q'}.
 \end{align}
Both integrals above are finite since $n<Q<p$ implies $nQ/p<n$ and $(n-1)Q'<n$. By a simple computation
 \be\label{Estim:Lrr:1}
  \Lambda_r(x)
   \leq
    C(n,p)\big(I_k[u;D]\big)^{1/p}
    \Big(\int_0^rt^{n-1-nQ/p}Y_{k+1}^{-Q/p}(t/D)\dd t\Big)^{1/Q}r^{n/Q'-n+1},
 \ee
for any $D\geq\max\{b''',\tilde{B}\}\diam(\Om)$. Lemma \ref{Lemma:X} for $\al=n-1-nQ/p$ and $\beta=Q/p$ ensures the existence of constants $\eta\geq0$ and $c>0$ both depending only on $n,p,Q$ (and thus only on $n,p)$, such that
 \[
  \int_0^rt^{n-1-nQ/p}Y_{k+1}^{-Q/p}(t/D)\dd t
   \leq c
    r^{n-nQ/p}Y_{k+1}^{-Q/p}(r/D),
 \]
for any $D\geq e^\eta\diam(\Om)$. Thus \eqref{Estim:Lrr:1} becomes
 \be\label{Estim:Lrr}
  \Lambda_r(x)
   \leq C_4(n,p)
    r^{1-n/p}Y_{k+1}^{-1/p}(r/D)\big(I_k[u;D]\big)^{1/p},
 \ee
for any $D\geq B\diam(\Om)$, where $B:=\max\big\{\max\{b''',\tilde{B}\},\eta\big\}$. Altogether, \eqref{Estim:Krr} and \eqref{Estim:Lrr} when inserted in \eqref{Estim:Krr:Lrr} give
 \[
  J(x)
   \leq C(n,p)
    r^{1-n/p}Y_{k+1}^{-1/p}(r/D)\big(I_k[u;D]\big)^{1/p},
 \]
for any $D\geq B\diam(\Om)$. The proof of \eqref{k:improvedHardy:Morrey} follows since the same estimate holds true for $J(y)$.

\medskip

To show the exponent $1/p$ on $X_{k+1}$ cannot be decreased, assume in contrary there exists $\e\in(0,1]$ such that for all $u\in\test(\Om\setminus\{0\})$ we have
 \be\nonumber
  \big(I_k[u;D]\big)^{1/p}
   \geq
    c\sup_{\substack{x,y\in\Om \\ x\neq y}}\Big\{\frac{|u(x)-u(y)|}{|x-y|^{1-n/p}}Y_{k+1}^{(1-\e)/p}\Big(\frac{|x-y|}{D}\Big)\Big\},
 \ee
for some constants $c>0$ and $D\geq\diam(\Om)$. Choosing $y=0$ we obtain
 \be\nonumber
  \big(I_k[u;D]\big)^{1/p}
   \geq
    c\frac{|u(x)|}{|x|^{1-n/p}}Y_{k+1}^{(1-\e)/p}(|x|/D),
     ~~~~\forall x\in\Om\setminus\{0\}.
 \ee
This readily implies that
 \be\label{violation:arg}
  \int_{\Om}\frac{|u(x)|^p}{|x|^p}Y_{k+1}^{2-\e/2}(|x|/D)\dd x
   \leq
    c^{-p}I_k[u;D]\int_{\Om}|x|^{-n}X_1^{1+\e/2}(|x|/D)\dd x.
 \ee
Clearly, since $\e>0$ the integral on the right is a finite constant depending only on $n,p,\e$ and $\Om$. Thus we have violated the optimality of the exponent $2$ of the remainder term (\ref{Hardy-improved:L^2}). \qed\\\\
\noindent {\textbf{Acknowledgements}} The first author was partially supported by Fondecyt grant 3140567 and by Millenium Nucleus CAPDE NC130017.

\author{\textsc{Konstantinos T. Gkikas}\\ Centro de Modelamiento Matem\'{a}tico (UMI 2807 CNRS)\\
 Universidad de Chile, Casilla 170 Correo 3, Santiago, Chile;\\ \texttt{kugkikas@gmail.com}\\
 \and\\
   \indent \textsc{Georgios Psaradakis}\\ Institut f\"{u}r Mathematik (Lehrst\"{u}hl  f\"{u}r Mathematik IV)\\
   Universit\"{a}t Mannheim, A5, Mannheim 68131, Deutschland;\\ \texttt{psaradakis@uni-mannheim.de}}
\end{document}